# A Two-Stage Stochastic Optimization Model for the Equitable Deployment of Fixed and Mobile Electric Vehicle Charging Stations

Hamid Najafzad [a], Moddassir Khan Nayeem [a], Fuhad Ahmed Opu [a], Omar Abbaas [a*], and Gabriel Nicolosi [b]

[a] Department of Mechanical, Aerospace, and Industrial Engineering, The University of Texas at San Antonio, San Antonio, TX 78249, USA

[b] Department of Engineering Management and Systems Engineering, Missouri University of Science and Technology, Rolla, MO

## Abstract

A major barrier to wide adoption of Electric Vehicles (EVs) is the absence of reliable and equitable charging infrastructure. Poorly located charging stations create coverage gaps and slow down EV adoption, especially in underserved communities. This paper proposes a two-stage stochastic mixed-integer programming model for the optimal deployment of Fixed and Mobile Charging Stations (FCSs and MCSs) across multiple zones and periods. Initially, a finite dominating set of candidate locations is identified using the Edge Scanning Algorithm for a Single Refueling Station (ESS), an exact, continuous-location method. We modify the ESS algorithm to incorporate existing public charging stations, thereby avoiding redundant coverage. In the first stage of our model, FCSs are allocated based on long-term traffic patterns, budgetary constraints, and socioeconomic factors to ensure stable baseline coverage. The second stage dynamically assigns MCSs in response to short-term demand fluctuations and uncertainties, aiming to minimize relocation costs while maximizing coverage. We use a scenario-based framework to capture demand variability. Numerical experiments on realistic networks demonstrate the model's capacity to enhance system resilience and reduce unmet demand. These findings offer practical insights for planners and policymakers seeking to develop accessible and demand-responsive EV charging infrastructure.

(*) Corresponding author, e-mail address: omar.abbaas@utsa.edu.

# 1. Introduction

The transportation sector is a major contributor to global carbon emissions and air pollution. According to the 2022 Inventory of U.S. Greenhouse Gas Emissions and Sinks, transportation accounted for nearly 28% of total U.S. greenhouse gas emissions, the largest share among all sectors [1-3]. In response, many countries have introduced stringent emission-reduction policies. Aligned with the Paris Agreement goal of limiting global warming to well below 2 °C, these policies have accelerated the global transition toward Electric Vehicles (EVs) as a key strategy for decarbonizing transportation [2, 4, 5].

In recent years, EV adoption has surged, with about 14 million electric cars sold globally in 2023, increasing the total EV fleet to nearly 40 million and accounting for an 18 percent share of all new car sales [6]. This growth has been most pronounced in China, Europe, and the United States, which together account for approximately 90% of global EV sales. Despite this expansion, EV penetration remains concentrated in these advanced markets, emphasizing the need for broader international and regional efforts to improve accessibility and supporting infrastructure. Widespread electrification of road transport is expected to prevent about two gigatons of $CO_2$-equivalent emissions by 2035 under current policies, even after including additional emissions from electricity generation [6]. In 2022, electric passenger vehicles in China accounted for about 35% of the global reduction in road transport emissions, showing the significant climate impact of early, large-scale EV adoption [6]. Despite these promising outcomes, the success of this transition depends on the availability of reliable and accessible charging infrastructure [7]. Early adopters have primarily relied on home-based charging, but the mass-market expansion of EVs requires a robust public charging network to ensure convenient access for all drivers.

Public charging capacity is projected to increase nearly sixfold by 2035, supporting large-scale EV adoption beyond major urban centers where most charging has historically been concentrated. According to the International Energy Agency (IEA), public chargers will increase from about 4 million in 2023 to roughly 25 million by 2035 [6]. Over the same period, private chargers will grow from about 36 million to more than 430 million [6]. The expansion of high-power public chargers is crucial to address range anxiety, which is the driver's concern about running out of charge before reaching a charging station. This infrastructure is essential for enabling long-distance travel and ensuring that EV adoption remains a practical option for all drivers [6]. Although public charging stations are expected to grow substantially, their current distribution remains highly uneven, concentrated in metropolitan centers and along major transport corridors, with limited availability in rural and low-density areas [8].

This imbalance is clear in real-world systems. Figure 1 shows the distribution of high-power public EV charging stations in Bexar County, Texas, where chargers mainly cluster in central and northern areas, while southern and outer regions remain underserved. This example demonstrates how uneven infrastructure deployment can create access gaps and worsen spatial inequalities, underscoring the need for strategic, data-driven planning.

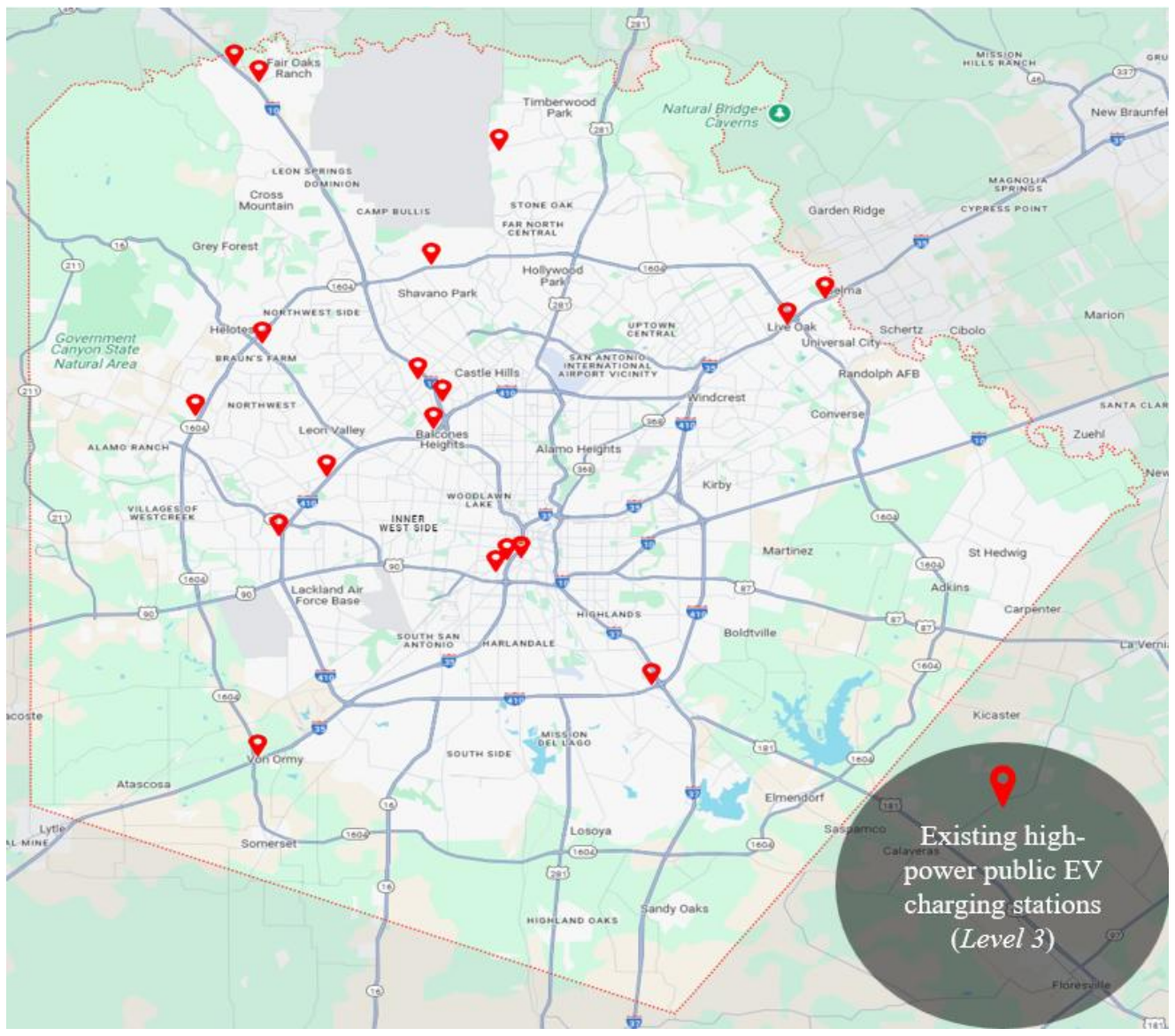

Figure 1: Current distribution of EV charging stations in Bexar County, Texas, USA

Expanding EV charging networks requires strategic planning by policymakers and industry leaders to determine the optimal locations, scale, and extent of deployment at the national, regional, and city levels [9]. Furthermore, proper management of battery charging resources and ongoing investment in transportation infrastructure are essential to handle increasing EV charging demand and ensure reliable system performance [10].

The placement of EV charging stations is crucial for the effectiveness of electric mobility systems [11]. Many existing charging networks are developed based on short-term profit market strategies, often prioritizing areas with high EV ownership while neglecting regions with lower adoption rates [12]. This uneven distribution leaves underserved areas with lower-income populations lacking adequate access to charging stations, thereby exacerbating socioeconomic inequities and slowing the transition to clean transportation [12]. These accessibility gaps emphasize the need for strategic, data-driven infrastructure planning. Ensuring that EV infrastructure is socially equitable is crucial, as underserved communities often face limited access, longer travel distances, and slower adoption rates. Incorporating equity into infrastructure design is therefore essential to achieving a fair and inclusive clean-mobility transition. Additionally, traditional location-allocation methods usually treat infrastructure as static and do not account for the potential integration of mobile charging technologies that can respond to short-term fluctuations [13]. As a result, while current methods maximize utilization in dense regions, they often leave considerable gaps in accessibility elsewhere.

To address these challenges, this study presents a two-stage stochastic optimization framework for the coordinated deployment of Fixed and Mobile Charging Stations (FCSs and MCSs) in an urban transportation network. The framework leverages the complementary roles of fixed and mobile infrastructure. FCSs are expensive to install and require sufficient grid capacity, making them suitable for locations with stable high charging demand. If underutilized, they can cause financial losses. MCSs offer a flexible and lower cost option. They can cover temporary gaps, respond to uncertain demand, and support market testing [14]. Additionally, MCSs can serve as an interim solution until permanent grid upgrades can support installing FCS. Building on this hybrid strategy, the framework incorporates socioeconomic factors to ensure equitable access to EV charging infrastructure across communities with different demographic and economic characteristics. Equity is measured through zone-specific weighting factors based on four key socioeconomic factors: Median household income, population density, public transit accessibility, and EV ownership levels. These weights are determined using the Best-Worst Method (BWM) [15], which systematically integrates expert judgment to assess the relative importance of each factor. The resulting equity index allows us to assign higher weights to underserved zones, minimizing spatial disparities in charging accessibility.

To improve candidate site selection for FCSs, a Modified ESS is proposed as a preprocessing step. Adapted from the Edge Scanning algorithm for a Single refueling station (ESS) introduced by Abbaas and Ventura [16], the Modified ESS extends the continuous location approach by incorporating existing public charging stations into the analysis. This ensures that new locations complement rather than duplicate existing infrastructure by identifying unique candidate points that address previously uncovered origin-destination pairs. These refined candidate locations are then used in the first-stage optimization to determine the optimal number and placement of FCSs. This stage focuses on maximizing flow coverage, equity, and accessibility across zones, while adhering to budget constraints. The second stage dynamically allocates Modified ESS across multiple time periods to handle the variations in traffic flow and charging demand. After demand scenarios are realized, MCSs are repositioned to high-demand areas to improve service availability, reduce relocation and operating costs, and ensure fair coverage over time. This adaptive stage enhances the network's flexibility and resilience, ensuring resources are used efficiently under uncertain, changing conditions.

To the best of the authors' knowledge, this paper is the first to introduce an integrated two-stage stochastic optimization framework that simultaneously addresses flow uncertainty, social equity, and the combined deployment of fixed and mobile EV charging stations. Unlike traditional deterministic or purely demand-based models, the proposed framework incorporates socioeconomic equity weights to promote equitable access to charging infrastructure across communities. It also presents the Modified ESS method, which combines existing charging facilities and identifies potential locations along the road network, reducing redundancy in station allocation. By integrating these elements, the framework provides a scalable, adaptive, and equity-focused decision-support tool for planners and policymakers to develop efficient, resilient, and socially inclusive EV charging networks.

The remainder of the paper is organized as follows. Section 2 reviews the related literature. Section 3 presents the problem statement. Section 4 describes the proposed methodology. Section 5 discusses computational experiments, covering data generation, experimental setup, and results analysis. Finally, Section 6 concludes the paper and demonstrates directions for future research.

## 2. Literature Review

The rapid adoption of EVs has increased demand for reliable and accessible charging infrastructure. This demand is influenced by factors like cost, charging time, station location, and grid impacts [3, 17]. Among these factors, the availability and spatial distribution of charging infrastructure are key in influencing EV adoption. Multiple studies demonstrate that infrastructure location is a more powerful predictor of adoption than financial incentives alone [18, 19]. As EV adoption grows, charging demand becomes closely tied to traffic patterns and power system constraints. Strategic EV infrastructure planning, therefore, focuses on optimal charger placement, service quality, and coordination with renewable energy sources to ensure grid stability and cost-effectiveness [20-22]. Understanding how users charge their vehicles is also crucial for effective planning. Empirical studies from Japan, the United Kingdom, and Ireland reveal notable differences in charging habits. Private users typically rely on home charging or fast-charging stations, whereas commercial users prefer workplace charging or stations that require minimal route changes [23, 24]. Moreover, peak charging periods emphasize the need for optimized routing, strategic station placement, pricing incentives, and social welfare-focused planning to balance grid demand and enhance network efficiency [25-28].

Despite these insights, most existing planning models remain static and deterministic. They usually assume fixed demand and consistent traffic conditions. The following subsections therefore review previous work on location-allocation models for EV charging infrastructure, demographic and socioeconomic factors influencing charging accessibility, mobile charging station systems, and stochastic optimization methods for handling flow uncertainty in transportation networks.

### *2.1. Location and Allocation Problem*

Location and allocation modeling plays a central role in EV charging network design, determining optimal station placement and capacity to balance accessibility, cost, and grid stability [3, 29-31]. Recent studies have integrated concepts like vehicle routing, detour reduction, and coordination with power distribution systems to improve service efficiency and ensure grid reliability [32, 33]. To address these challenges, different modeling paradigms such as flow-based, node-based, path-based, tour-based, and network equilibrium have been developed to capture the complex interactions between transportation demand, charging accessibility, and energy supply systems [3].

Flow-based models are particularly effective in identifying high-impact locations along major travel routes by capturing vehicle movement rather than static demand. One of the pioneering works in flow-based models is the flow-capturing location model (FCL), proposed by Hodgson [34]. The FCL model and its extensions [35-38] focus on high-traffic areas to improve coverage and reduce detours by including concepts such as vehicle range limits, multi-stop refueling, and spatial-temporal demand changes to boost station placement and charging network effectiveness. Building on this foundation, numerous studies across countries such as the United States, South Korea, Spain, and China have used flow-based frameworks to optimize station placement, lower operational costs, and improve accessibility [39-43]. However, these models are typically formulated on discrete networks, where candidate sites are restricted to existing nodes. Such discrete representations may lead to suboptimal placement because potential sites between nodes are ignored. To overcome this limitation, researchers introduced continuous network models that consider every point along road segments as potential station locations [44]. These models allow detours for refueling, increasing route coverage, and ensuring global optimality, a property that a discrete method cannot guarantee [45]. Further, Abbaas and Ventura [16] developed the ESS, a polynomial-time exact method that systematically identifies globally optimal refueling locations along network edges. Later refinements extended this framework to handle multiple stations and lexicographic objectives [46].

A node-based strategy efficiently allocates and optimizes network resources by organizing charging stations as interlinked nodes [3]. Each node symbolizes a single station or a group, managed centrally to oversee electricity flow, modify pricing, and distribute power according to real-time demand [3]. Many studies have employed node-based methods to enhance the locations of EV charging stations, considering factors such as population, demand, and travel time [47-49]. Other studies have further improved these models by integrating demand forecasting, traffic flow, and multi-period planning to boost scalability and long-term network efficiency [50-52]. Conversely, a path-based approach enhances EV charging infrastructure by placing stations along vehicle routes to provide timely access, reduce congestion, and ease range anxiety. Studies using path-based models have shown improved traffic flow management, reduced system costs, and enhanced network accessibility across various application scenarios [3, 53-55].

Tour-based approaches extend these models by planning routes that include necessary charging stops rather than treating stations as isolated points [3]. These models improve efficiency and user experience by minimizing stops and congestion while considering battery level, cost, and station availability [3]. Various studies have applied tour-based models to optimize station siting, sizing, and routing by minimizing deviations, lowering installation costs, and improving accessibility to ensure efficient and reliable charging network performance [56-58]. Network equilibrium frameworks optimize EV charging infrastructure by balancing demand and supply while minimizing costs and maximizing network efficiency [3]. These models integrate traveler behavior, energy constraints, and flow-dependent factors to reduce trip times, alleviate congestion, and enhance social welfare [59-63]. Other studies have extended this framework to coordinate transportation and power systems, enhancing cost efficiency, service quality, and infrastructure accessibility [64-67]. While these models significantly improve EV infrastructure planning, most overlook the impact of existing charging networks and fail to account for stochastic variations in traffic flow and demand, key gaps that the current paper seeks to address.

### *2.2. Demographic and Socioeconomic Influences on EV Infrastructure*

Demographic and socioeconomic factors, such as income, race, education, and housing, significantly shape the distribution of EV charging infrastructure, often resulting in unequal access across communities [68-72]. Studies show that charging stations are concentrated in affluent, white-majority neighborhoods, while low-income, minority, and multi-unit housing areas face limited access, particularly in cities like New York, Chicago, and across California [68-74]. Economic activity, political affiliation, and proximity to highways further influence charger placement, favoring wealthier and urban areas [75]. Moreover, high installation costs, grid limitations, and limited home-charging options exacerbate disparities, especially in rural and disadvantaged areas [76-78]. To mitigate such inequalities, recent studies have suggested including equity factors, such as income level, population density, or transit access, in charging station location models [79, 80]. However, most current methods rely on post-hoc fairness assessments or address equity qualitatively instead of incorporating it into optimization models. This shows the need for quantitative, equity-centered optimization approaches that explicitly include socioeconomic factors in infrastructure planning decisions.

### *2.3. Mobile Charging Station*

Addressing the increasing and unpredictable energy needs of EVs poses a significant challenge in the shift toward an electric transportation system [81]. Charging needs fluctuate with local congestion, travel behavior, and seasonal patterns [82, 83], requiring flexible and adaptive infrastructure. These complexities require flexible, time-sensitive charging solutions. MCSs have become a cost-effective complement to fixed infrastructure, especially in underserved or high-demand areas [84-88]. Their mobility allows for quick deployment, relocation, and grid integration, helping to ease congestion at fixed stations, lessen range anxiety, and promote equitable access [89-93]. By dynamically responding to spatial and temporal demand variations, MCSs optimize network coverage, minimize wait times, and improve user satisfaction [94-96]. Furthermore, they promote grid stability and sustainability by coordinating with renewable energy sources like solar and wind power [97, 98]. Recent advancements, including bilevel and Lyapunov-based optimization frameworks, improve dynamic routing, scheduling, and battery management to enhance operational efficiency, energy distribution, and profitability in real-time environments [99-105]. Despite these advances, existing studies often consider MCS deployment separately from fixed infrastructure and neglect stochastic flow variations. This work addresses that gap by developing an integrated two-stage stochastic framework that jointly allocates FCSs and MCSs to enhance flexibility, resilience, and equitable network coverage.

### *2.4. Stochastic Programming for EV Infrastructure*

Uncertainty makes deterministic models less effective and more costly in unpredictable environments [105]. To address this, stochastic programming uses probability distributions to capture variability and enhance solution robustness [106]. Various approaches, such as robust, probabilistic, fuzzy, and Bayesian optimization, have been successfully used in fields like logistics, supply chain, energy, and production planning [107-111]. In the EV sector, stochastic models effectively handle fluctuations in energy procurement, dynamic pricing, and mobility demand, enhancing reliability and decision-making efficiency [112, 113].

Two-stage stochastic programming is commonly used to handle uncertainty by dividing decisions into first-stage (pre-uncertainty) and second-stage (adaptive) actions based on realized outcomes [114, 115]. This scenario-based framework improves flexibility but has seen limited use in flow refueling, where demand and supply conditions are uncertain [116]. Recent studies have incorporated demand uncertainty into flow-interception and dynamic service models for EV charging and battery-swapping stations [117, 118], optimized power distribution under uncertain demand [119], and improved charging station placement for demand variability [120]. Two-stage approaches have also been used to integrate demand learning and dynamic deployment for refueling stations [121], while other works address uncertainties in energy transactions [122], driving range [123, 124], battery status, power consumption [125, 126], and stochastic travel behavior [127-129]. Collectively, these studies confirm that stochastic optimization, especially two-stage frameworks, improves the adaptability and robustness of EV charging infrastructure planning under uncertainty.

*2.5. Research Gaps*

The reviewed literature shows significant progress in charging infrastructure planning but also reveals essential limitations. Existing studies often ignore how the current charging infrastructure is integrated, treat location and allocation decisions separately, and depend on deterministic assumptions that do not account for flow uncertainty and fairness issues. To address these gaps, this study presents two main contributions. First, we improve the ESS by adding existing public charging stations into the network to avoid redundancy and focus on underserved areas. This change increases placement efficiency, lowers installation costs, and boosts computational performance by targeting high-impact candidate sites. This framework expands the Continuous Deviation-Flow Refueling Station Location (CDFRSL) [16] problem by incorporating multi-zone, multi-period, and equity-focused planning within a unified stochastic model. The resulting approach offers a scalable, data-driven, and socially inclusive solution for optimizing real-world EV charging infrastructure. Second, to manage allocation under uncertainty, we develop a two-stage stochastic mixed-integer programming (MIP) model that jointly optimizes the deployment of FCSs and MCSs across multiple zones and time periods. The first stage ensures equitable access by allocating FCSs based on average flow and socioeconomic weights. In contrast, the second stage dynamically repositions MCSs under different demand scenarios to minimize relocation costs and improve network coverage.

## 3. Problem Statement

In this paper, we introduce an optimization framework to identify optimal locations and allocation strategies for FCSs and MCSs within an urban transportation network. We consider a connected, undirected traffic network denoted by $G(V, E)$, where $V$ represents a set of $n$ vertices (interchanges) and $E$ represents a set of $e$ edges (road segments). An edge connecting vertices $v_a$,$v_b \in V$ is denoted by $(v_a, v_b)$. This edge is undirected, meaning that vehicles can travel from $v_a$ to $v_b$ and vice versa, and the distance is the same in both ways. To prevent repetition, edges are included only when $a < b$.

In this network, vehicles travel from an origin vertex $v_i \in V$ to a destination vertex $v_j \in V$ using the shortest refuellable path, a path that includes at least one available charging station within the vehicle's travel range. After reaching $v_j$, vehicles complete a round trip by returning to $v_i$ via the same path in the reverse direction. Let $q(v_i, v_j)$ represent an Origin-Destination (O-D) pair of vertices, where $v_i$,$v_j \in V$ and $i < j$. The average traffic flow for $q(v_i, v_j)$, is denoted by $f(v_i, v_j)$. This flow, measured in roundtrips per time unit, reflects the traffic in both directions between the O-D pair. $G(V, E)$ is a general network; hence there could be several paths between any two points in the network. The shortest path between points $x_1$ and $x_2$ is denoted by $P(x_1, x_2)$, with its length referred to as $d(x_1, x_2)$. A closed line segment between two points, $x_1$ and $x_2$, on the same edge is referred to as $C(x_1, x_2)$, with the interior points denoted by $int(C(x_1, x_2))$, and the length of the segment is $l(x_1, x_2)$.

This paper presents a hybrid optimization framework designed for the fair and efficient placement of FCSs and MCSs across a multi-zone, multi-period traffic network. It accounts for spatial and temporal dynamics to address varying EV demand while balancing accessibility, efficiency, and equity. Fixed stations provide stability in high-demand areas, whereas mobile stations offer flexibility to accommodate short-term demand spikes. By leveraging existing infrastructure and a continuous location method to identify candidate endpoints, the model reduces redundancy, enhances coverage, and improves computational efficiency, thereby maximizing network-wide flow coverage and ensuring fair access to charging in underserved regions.

Charging needs vary across zones and time periods due to changing traffic patterns and unpredictable user behavior. To address this, the model uses a scenario-based stochastic approach that explicitly accounts for demand and resource-availability uncertainty. This framework optimizes both FCS placement and MCS deployment under various possible conditions, ensuring adaptability, resilience, and cost-effectiveness. By combining stochastic elements with financial and social goals, the proposed methodology supports robust infrastructure planning that remains effective amid real-world variability. Assumptions for this problem include the following:

I. A charging station (FCS or MCS) has no capacity constraints and can accommodate vehicles traveling in both directions.

II. All vehicles have the same travel range, referred to as $R$.

III. Vehicles travel in round trips between O-D pairs. A complete round trip consists of an original trip from the origin vertex to the destination vertex and a return trip from the destination vertex to the origin vertex using the same path in the reverse direction.

IV. Both original and return trips start with a charge enough to drive at least $\alpha R$, where $0 \leq \alpha \leq 1$. Moreover, a vehicle is allowed at most one refueling operation in each direction.

V. The distance between any two points, $x_1, x_2$, is symmetric; therefore, $d(x_1, x_2) = d(x_2, x_1)$.

VI. Drivers can deviate from their shortest path $P(v_i, v_j)$ to access a charging station, without exceeding the vehicle's range.

The first assumption ensures that each charging station, FCS or MCS, can provide service to all vehicles traveling in both directions. The second assumption states that all EVs share the same travel range $R$, making charging needs directly proportional to the distance traveled. The third assumption defines O-D flows as complete round trips, which reflects typical travel behavior in many transportation settings, such as commuting from home to work and returning. The fourth assumption requires that each EV begins its trip with sufficient energy to travel at least $\alpha R$, thereby maintaining a minimum state-of-charge that ensures reliable operation in both travel directions. This assumption implies that the length of any edge $(v_i, v_j) \in E$ must be less than or equal to $2\alpha R$, ensuring that vehicles can start their trips from any of its vertices and access a charging station without exceeding their range. Moreover, for most intracity applications and under reasonable full-charge driving ranges, it is unlikely that a vehicle would require multiple recharging events during a single one-way trip. The fifth assumption states that distance between any two points in the network is symmetric. The sixth assumption reflects realistic driver behavior, allowing EVs to take detours from the shortest path $P(v_i, v_j)$ to reach a nearby charging station without exceeding their driving range.

## 4. Methodology

This section describes the proposed methodology. Subsection 4.1 explains the proposed Modified ESS algorithm, which considers existing charging stations and identifies candidate locations for charging stations. Subsection 4.2 employs BWM to determine zone-specific socioeconomic weights that represent equity considerations within the model. Subsection 4.3 introduces the two-stage stochastic mathematical model that optimizes the locations and allocation of FCSs and MCSs across zones and time periods.

### *4.1. Modified Edge Scanning Algorithm for a Single Refueling Station*

Abbaas and Ventura [16] introduced the ESS to solve the CDFRSL problem on a general network. The algorithm systematically scans each network edge and, for every O-D pair, identifies refueling segments, defined as continuous portions of the edge where a refueling station can feasibly serve the flow between that O-D pair given vehicle range and deviation limits. Rather than evaluating all points along an edge, ESS establishes that the overall set of endpoints of these refueling segments always contain an optimal solution that maximizes flow coverage. As a result, ESS reduces the continuous problem to a finite set of candidates while preserving global optimality. The set of O-D pairs, referred to as $Q^o$, is defined as:

$$Q^o = \{q(v_i, v_j) \mid d(v_i, v_j) \leq 2\alpha R, f(v_i, v_j) > 0, i < j, v_i, v_j \in V\}. \quad (1)$$

An O-D pair $q(v_i, v_j)$ is considered covered by a charging station located at point $x$ if this point falls within $\alpha R$ from $v_i$ and $v_j$. The set of O-D pairs covered by a charging station located at point $x$, denoted by $S^o(x)$, is defined as follows:

$$S^o(x) = \{q(v_i, v_j) \mid d(v_i, x) \leq \alpha R, d(x, v_j) \leq \alpha R, i < j, q(v_i, v_j) \in Q^o\}. \quad (2)$$

The total traffic flow covered by a charging station at $x$ is denoted by $F^o(x)$:

$$F^o(x) = \sum_{q(v_i, v_j) \in S^o(x)} f(v_i, v_j). \quad (3)$$

In this research, the ESS algorithm is extended to improve computational efficiency and reduce redundant coverage. The Modified ESS incorporates existing charging infrastructure into the optimization process by identifying areas already served and focusing new station placement on underserved regions. The following steps outline the process for identifying the final set of candidates locations within the network.

Step 1. Identify the Domain of Existing Stations: Let $U$ denote the set of existing charging stations in the network. Each station $u \in U$ located at point $x_u$ in the network covers a set of O-D pairs, $S(x_u)$. The overall coverage of existing stations, representing the combined domain $Dom(U)$, is given by:

$$Dom(U) = \cup_{u \in U} S^o(x_u). \tag{4}$$

This includes every O-D pair that can be served by at least one existing station.

Step 2. Exclude Covered Pairs: After determining $Dom(U)$, the already-served O-D pairs are excluded from further analysis to concentrate computational effort on the remaining uncovered demand. The set of uncovered O-D pairs Q is defined as:

$$Q = Q^o \backslash Dom(U). \tag{5}$$

Step 3. Apply the ESS Algorithm for Uncovered O-D Pairs: We now apply the ESS algorithm introduced in [16] to the set of uncovered O-D pairs. For completeness, a brief overview of the ESS algorithm is provided below.

The ESS algorithm classifies vehicle trips to a charging station into two types based on how the vehicle traverses the edge containing the station. In Type 1 trips, vehicles reach a charging station on edge $(v_a, v_b) \in E$ by going through the entire edge. In this type, a vehicles enter the edge either using $v_a$ and exit using $v_b$ (Figure 2) or enter using $v_b$ and exit using $v_a$ (Figure 3). The maximum possible travel distance along edge $(v_a, v_b) \in E$ before refueling depends on its remaining driving range upon entering the edge and the edge length. This distance is denoted by $\gamma(v_k; v_r, v_c)$, where $v_k$ is the origin vertex in the original trip or destination vertex in the return trip, $v_r$ is the entry vertex to the edge containing the refueling station, and $v_c$ is the other vertex of the same edge. $\gamma(v_k; v_r, v_c)$ is defined as follows:

$$\gamma(v_k; v_r, v_c) = \min\{\alpha R - d(v_k, v_r), l(v_a, v_b)\},\ k = i, j; r = a, b; c \neq r. \tag{6}$$

The feasible refueling set for vehicles undertaking Type 1 trips is then expressed as:

$$\begin{aligned} &RS_1^{(r)}(v_i, v_j; v_a, v_b) = \\ &\quad \{x \in (v_a, v_b) \,|\, l(v_r, x) \leq \gamma(v_i; v_r, v_c)\ and\ l(v_c, x) \leq \gamma(v_j; v_c, v_r)\}, \\ &\quad r = a, b; c = a, b; c \neq r. \end{aligned} \tag{7}$$

By combining both directional scenarios (a) and (b), the complete Type 1 refueling segment is obtained as:

$$RS_1(v_i, v_j;\ v_a, v_b) = RS_1^{(a)}(v_i, v_j; v_a, v_b) \cup RS_1^{(b)}(v_i, v_j; v_a, v_b). \qquad (8)$$

As shown in Figure 4 and Figure 5, Type 1 refueling set may include up to four endpoints.

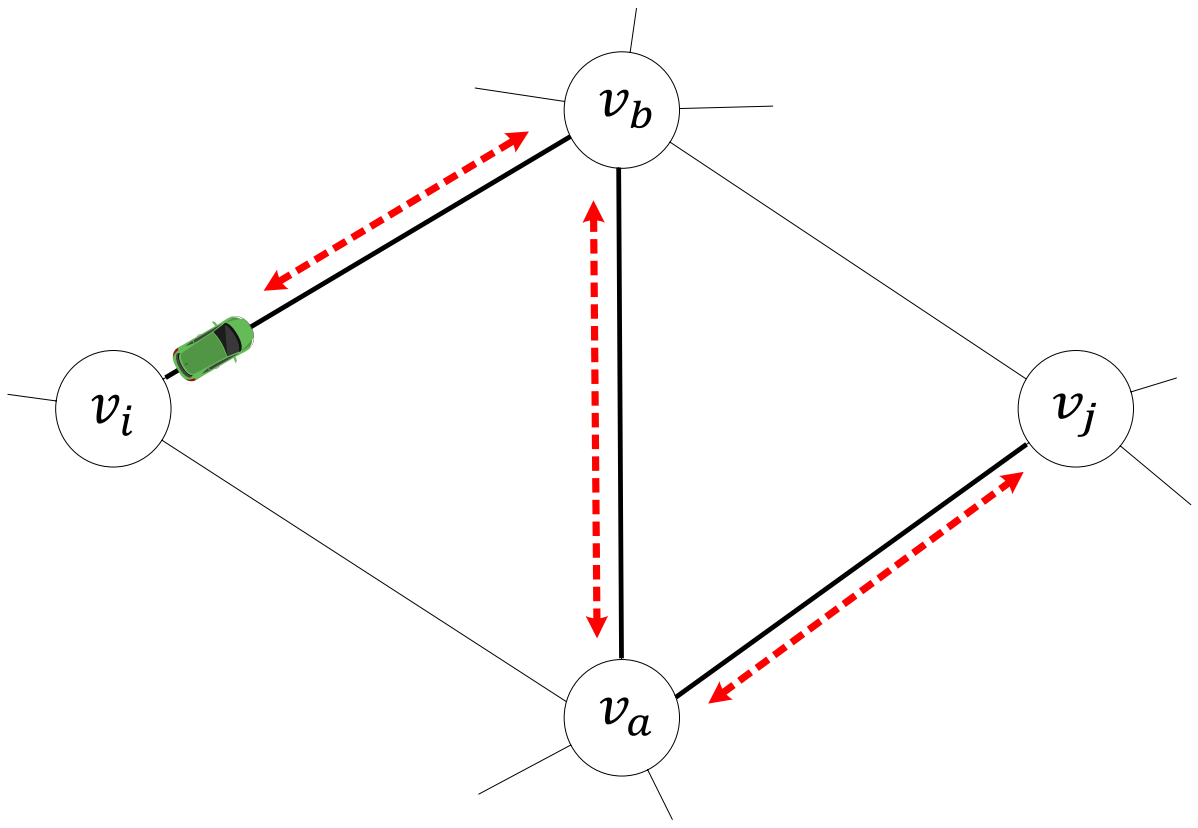

Figure 2: Type 1 scenario (a): Flow of original trips from $v_i$ to $v_j$ via $v_a$ and $v_b$

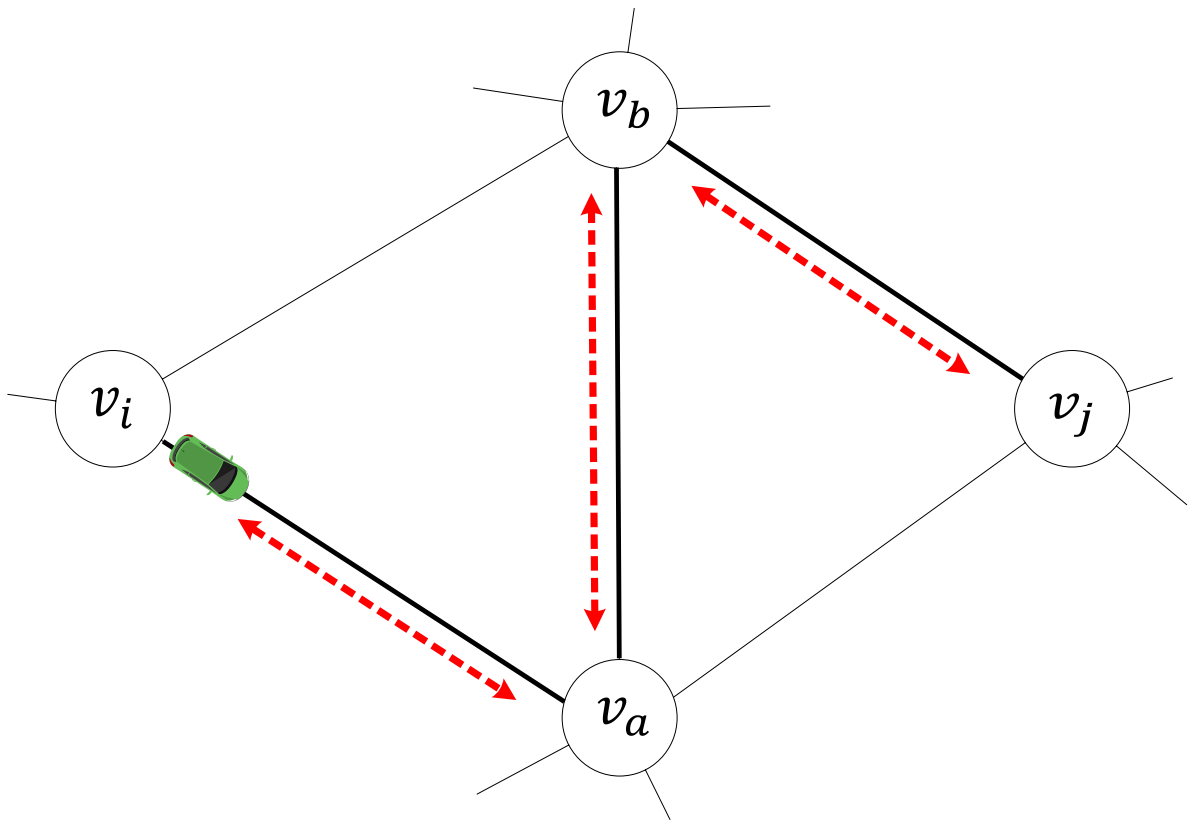

Figure 3: Type 1 scenario (b): Flow of original trips from $v_i$ to $v_j$ via $v_b$ and $v_a$

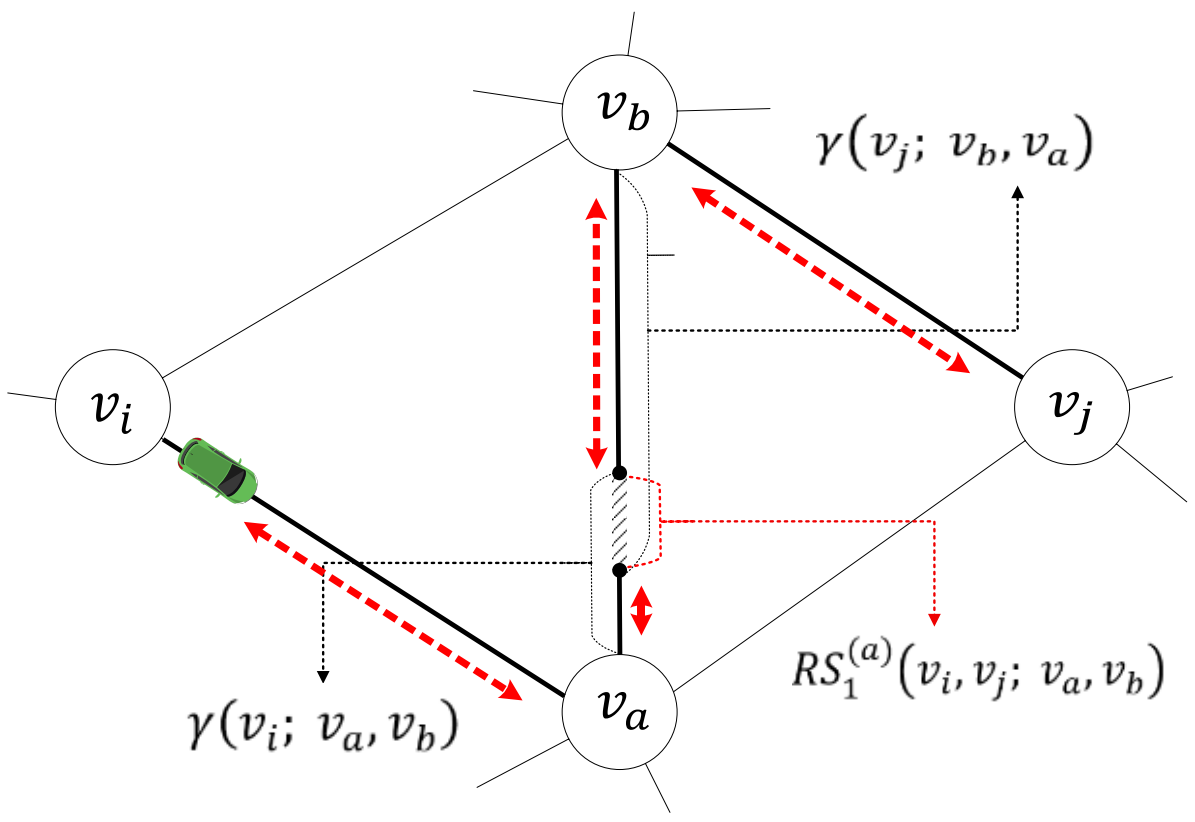

Figure 4: Refueling segment for scenario (a) of Type 1

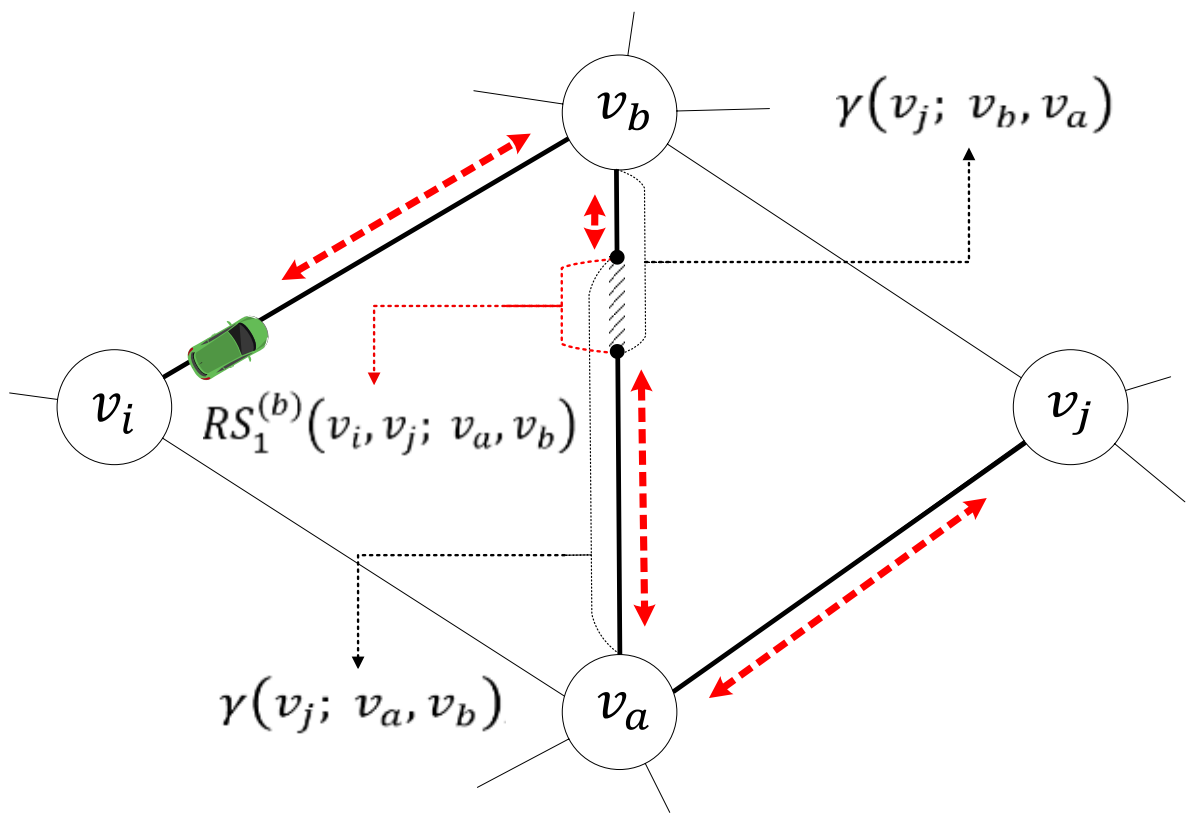

Figure 5: Refueling segment for scenario (b) of Type 1

In Type 2 trips, vehicles refuel along an edge $(v_a, v_b) \in E$. However, in this case, a vehicle enters the edge, travels to the charging station, makes a U-turn, and exits the edge using the same entry vertex. Vehicles may enter the edge using either $v_a$ (Figure 6) or $v_b$ (Figure 7). The maximum possible travel distance along edge $(v_a, v_b) \in E$ before refueling depends on its remaining driving range upon entering the edge and the edge length. This distance is denoted by $\delta^{(r)}(v_i, v_j; v_a, v_b)$, where $q(v_i, v_j) \in Q$, $(v_a, v_b) \in E$ and $v_r \in \{v_a, v_b\}$ is the entry vertex to the edge containing the refueling station. $\delta^{(r)}(v_i, v_j; v_a, v_b)$ is defined as follows:

$$\delta^{(r)}(v_i, v_j; v_a, v_b) = \min\{\alpha R - \max\{d(v_i, v_r), d(v_j, v_r)\}, l(v_a, v_b)\}, \quad r = a, b, \qquad (9)$$

And the corresponding refueling segment is defined in Equation ((10):

$$RS_2^{(r)}(v_i, v_j; v_a, v_b) = \{x \in (v_a, v_b) | l(v_r, x) \leq \delta^{(r)}(v_i, v_j; v_a, v_b)\}, \quad r = a, b. \tag{10}$$

The overall refueling set for both scenarios is given by Equation ((11).

$$RS_2(v_i, v_j; v_a, v_b) = RS_2^{(a)}(v_i, v_j; v_a, v_b) \cup RS_2^{(b)}(v_i, v_j; v_a, v_b). \tag{11}$$

Figure 8 and Figure 9 illustrate Type 2 refueling segments for both scenarios, showing how vehicles diverge from and return to the same vertex while completing their trip.

After identifying the refueling sets for both trip types, the overall refueling set for an uncovered O-D pair $q(v_i, v_j) \in Q$ on edge $(v_a, v_b) \in E$ is:

$$RS(v_i, v_j; v_a, v_b) = RS_1(v_i, v_j; v_a, v_b) \cup RS_2(v_i, v_j; v_a, v_b). \tag{12}$$

Each combined set may contain up to eight endpoints, four from Type 1 trips and four from Type 2 trips, denoted as $w^k_{v_i, v_j; v_a, v_b}$ for $k = 1, \dots, 8$ as expressed in Equation ((13).

$$EP(v_i, v_j; v_a, v_b) = \left\{ w^k_{v_i, v_j; v_a, v_b} \middle| k = 1,2 \dots, 8 \right\}. \tag{13}$$

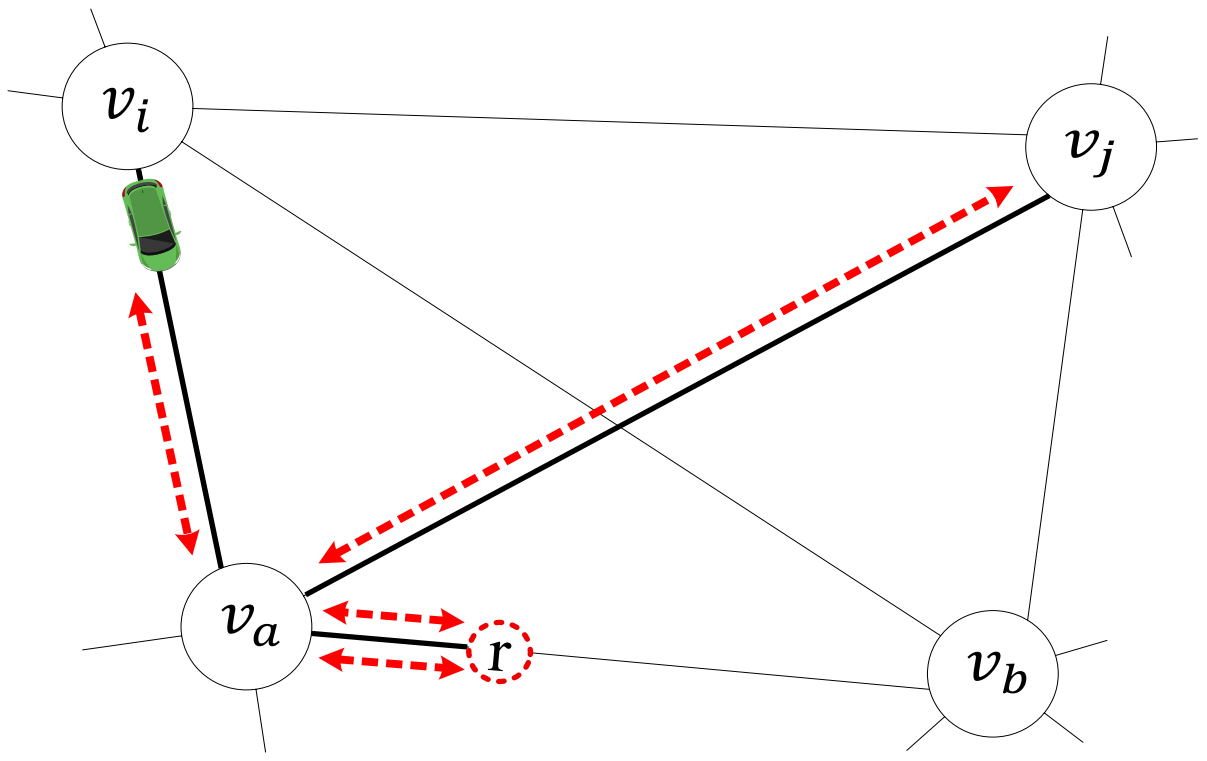

Figure 6: Type 2 scenario (a); Original Trips from $v_i$ to $v_j$: deviating from and returning to $v_a$

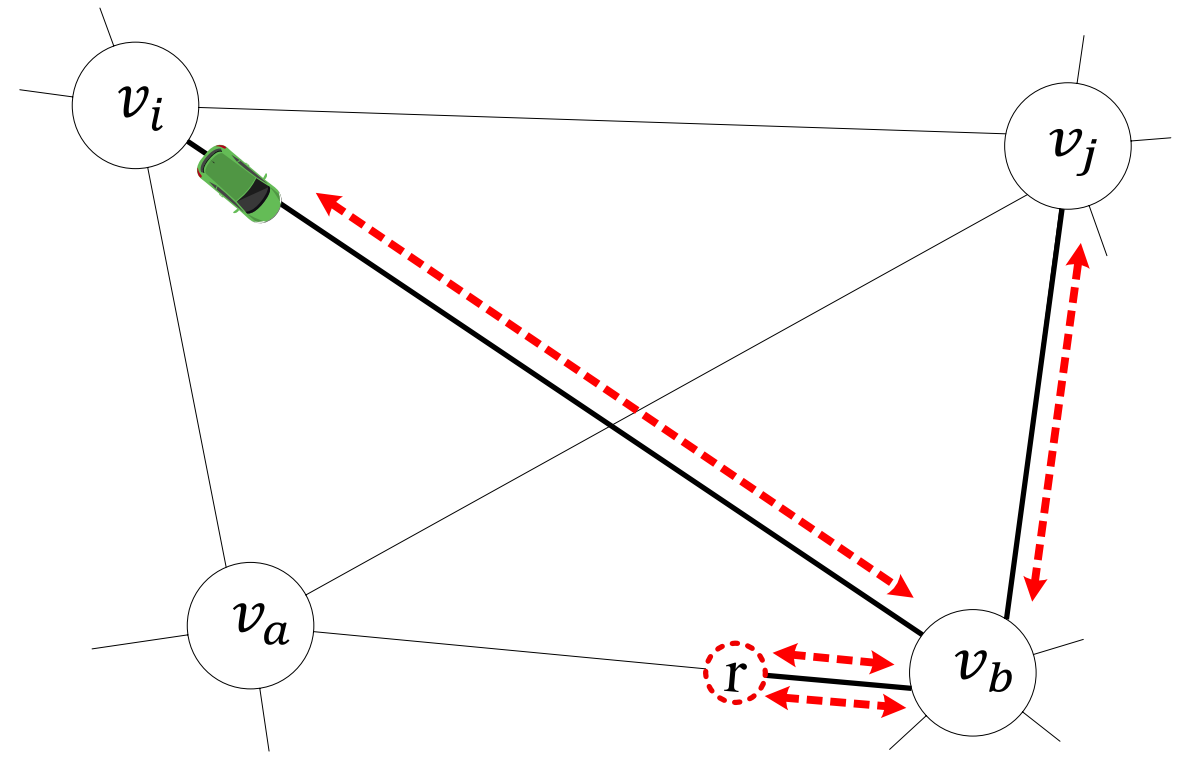

Figure 7: Type 2 scenario (b); Original Trips from $v_i$ to $v_j$: deviating from and returning to $v_b$

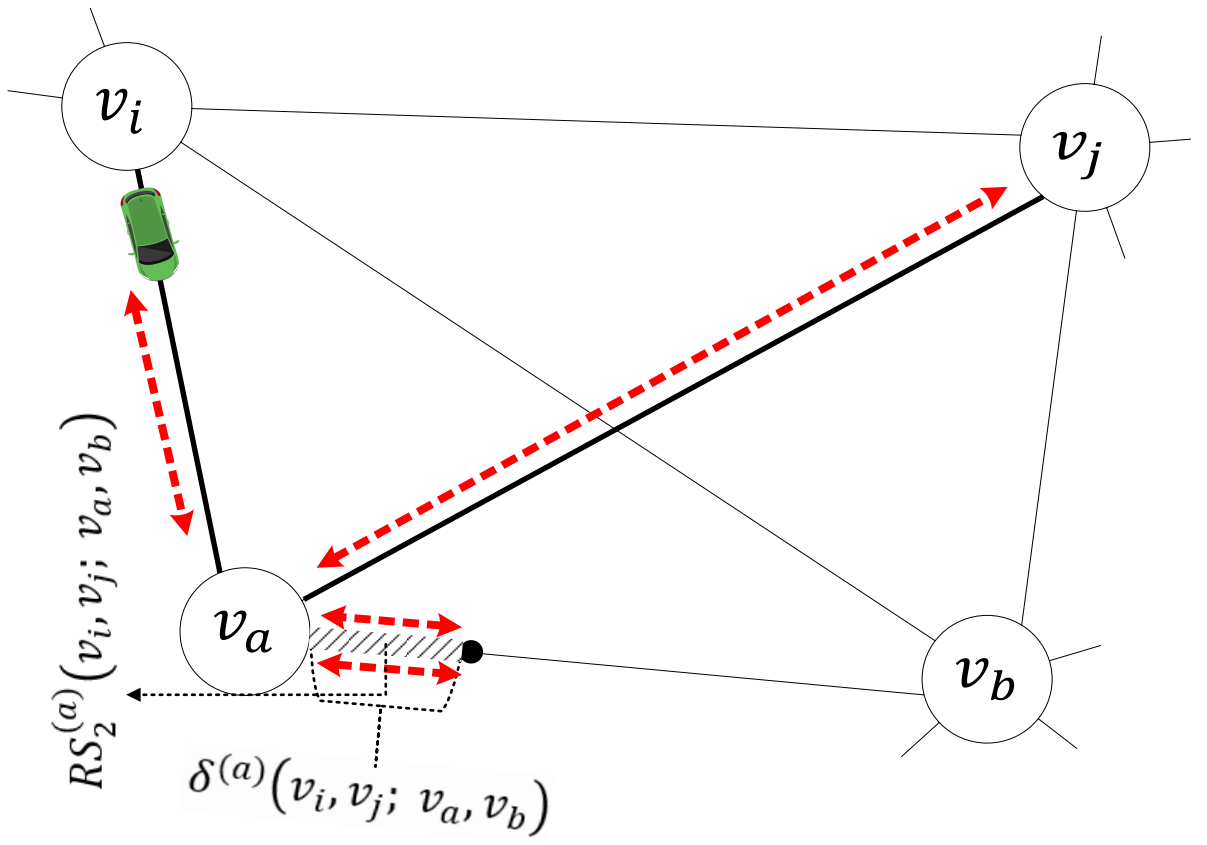

Figure 8: Refueling segment for scenario (a) of Type 2

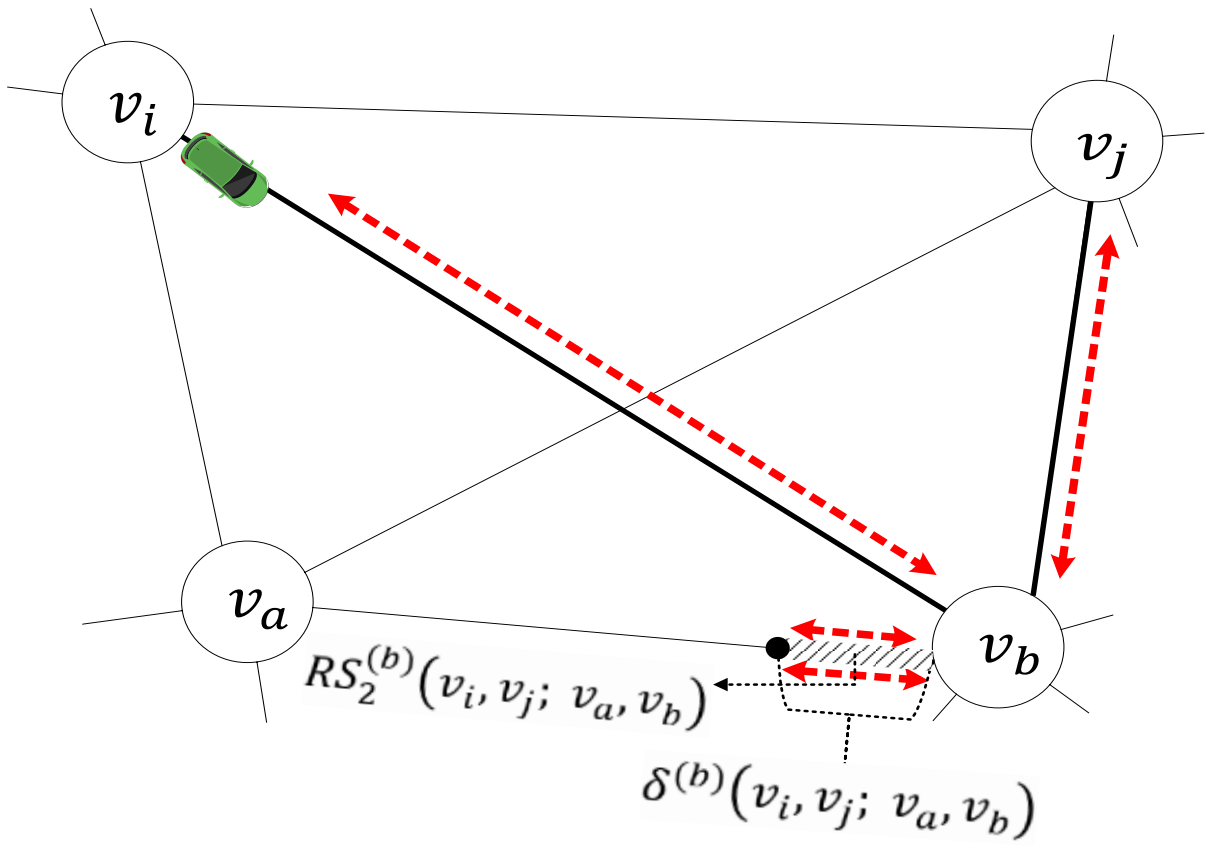

Figure 9: Refueling segment for scenario (b) of Type 2

As we consider all O-D pairs for potential refueling sets on edge $(v_a, v_b) \in E$, the set of endpoints on this edge, $EP(v_a, v_b)$, is defined as:

$$EP(v_a, v_b) = \cup_{q(v_i,v_j)\in Q} EP\left(v_i, v_j;\ v_a, v_b\right). \quad (14)$$

The overall set of endpoints for these uncovered O-D pairs is defined as follows:

$$EP = \cup_{(v_a,v_b)\in E} EP(v_a, v_b). \quad (15)$$

Refueling segments may overlap for several O-D pairs, this allows some endpoints to cover the flows between multiple O-D pairs. A charging station located at endpoint $w \in EP$, can cover a set of O-D pairs, $S(w)$, and a total flow, $F(w)$, defined as follows:

$$S(w) = \left\{q\left(v_i, v_j\right)\middle| d(v_i, w) \le \alpha R, d\left(w, v_j\right) \le \alpha R, i < j, q\left(v_i, v_j\right) \in Q\right\}, \quad (16)$$

$$F(w) = \textstyle\sum_{q(v_i,v_j)\in S(w)} f\left(v_i, v_j\right). \quad (17)$$

*4.2. Socioeconomic Weight Derivation Using the Best-Worst Method*

To integrate social equity into the optimization framework, we divide the network into a set of zones, denoted by $Z$, based on socioeconomic factors. The BWM [15] is then applied to derive zone-specific weights, which are subsequently embedded in the mathematical formulation presented in the next subsection. The weighting process incorporates four key factors: median household income, population density, access to public transportation, and EV ownership levels. This ensures that charging infrastructure deployment reflects both community needs and adoption potential within each zone. Areas with lower incomes and limited transit access stand to benefit most from EVs' lower operating costs and available government incentives. In addition, high population density and existing EV ownership signal strong urban demand, future growth potential, and a reduced likelihood of access to home charging. To achieve fair and data-driven planning, a weighting factor is assigned to each zone based on its socioeconomic profile. The BWM evaluates the significance of each factor by comparing the most and least important criteria, thereby reducing bias and ensuring consistency. It combines expert judgment with quantitative data to support equitable distribution of EV infrastructure across zones. Below are the summarized steps to apply BWM in this context.

Step 1. Identify the Best and Worst Criteria: Decision-makers identify the most and least important socioeconomic factors (e.g., median income as the most important and EV ownership as the least important).

Step 2. Best-to-Others Comparison: The most important (best) factor is rated on a scale from 1 to 9, where 1 means equal importance and 9 indicates extreme importance. This produces the Best-to-Others vector:

$$O_b = [o_{b1}, o_{b2}, o_{b3}, \dots, o_{bI}], \quad (18)$$

where $I$ represents the number of factors, while $o_{bi}$ indicates the relative importance of the best factor to the $i^{th}$ factor.

Step 3. Comparison Against the Worst Factor: All factors are then compared against the least important (worst) factor using the same scale, creating the Others-to-Worst vector:

$$O_w = [o_{1w}, o_{2w}, o_{3w}, \dots, o_{Iw}], \tag{19}$$

where $o_{iw}$ is the relative importance of the $i^{th}$ factor compared to the worst factor.

Step 4. Optimization Model: An optimization model is solved to determine the coefficients $\theta_i$ which represent the relative importance of factor $i \in \{1, \dots, I\}$. This model minimizes the maximum deviation between expert judgements and derived weight ratios as follows:

$$\text{Min } \varepsilon \tag{20}$$

Subject to:

$$\left|\frac{\theta_b}{\theta_i} - o_{bi}\right| \leq \varepsilon, \forall i \in \{1, \dots, I\}, \tag{21}$$

$$\left|\frac{\theta_i}{\theta_w} - o_{iw}\right| \leq \varepsilon, \forall i \in \{1, \dots, I\}, \tag{22}$$

$$\textstyle\sum_{i=1}^{I} \theta_i = 1, \theta_i \geq 0, \forall i \in \{1, \dots, I\}, \tag{23}$$

$$\varepsilon \geq 0. \tag{24}$$

The objective function ((20) minimizes $\varepsilon$, which represents the maximum deviation between the calculated weight ratios and expert-provided comparisons. Constraint sets ((21-(22) enforce consistency with the Best-to-Others and Others-to-Worst judgments, respectively. Constraint set ((23) normalizes the weights so they sum to one and remain nonnegative, ensuring valid proportional importance. Constraint ((24) ensures that the deviation variable $\varepsilon$ is nonnegative.

Step 5. Calculate Zone-Specific Weights: Let $\delta_i^z$ be the observed value of socioeconomic factor $i$ in zone $z$, and let $\delta_i^{min}$ and $\delta_i^{max}$ represent the minimum and maximum observed values of factor $i$ across all zones, respectively. Since the considered factors are measured on different scales, we normalize them as follows:

$$G_i^z = \frac{\delta_i^z - \delta_i^{min}}{\delta_i^{max} - \delta_i^{min}}, \tag{25}$$

where $G_i^z$ is the normalized value of the $i^{th}$ socioeconomic factor in zone $z$. This normalization maps all values to the interval [0,1], ensuring comparability across factors. To maintain consistent interpretation, all factors are scaled such that higher values indicate greater need or priority. For factors such as income and transit access, where higher values imply lower vulnerability, the scale is inverted using $1 - G_i^z$. This transformation ensures that all factors contribute meaningfully and consistently to the equity weighting process, aligning the data with planning objectives.

Next, the normalized factor values $G_i^z$ and their corresponding importance coefficients $\theta_i$ are used to calculate the composite socioeconomic weight for each zone $z \in Z$, denoted by $\mu_z$, as follows:

$$\mu_z = \sum_{i=1}^{I} \theta_i G_i^z, \forall z \in Z, \quad (26)$$

This procedure ensures that zones with greater socioeconomic disadvantage receive higher composite weights ($\mu_z$), thereby prioritizing them in the subsequent optimization process. Because our objective is defined overflows between O-D pairs, we also require a socioeconomic weight for each O-D pair $q(v_i, v_j)$, denoted by $\mu_{q(v_i,v_j)}$. However, it is possible for the origin and destination vertices to be in two different zones. Moreover, the path between the O-D pair could pass through multiple zones. Several mappings from zone weights to O-D weights are possible. For example, $\mu_{q(v_i,v_j)}$ could be computed as a distance-weighted average of the zone weights along the O-D path, or as the average of the origin- and destination-zone weights. In this study, we adopt a conservative prioritization rule and set $\mu_{q(v_i,v_j)}$ to the maximum of the origin and destination zone weights, so that O-D pairs involving at least one highly disadvantaged zone receive higher priority in the coverage objective. Let $V_z$ be the set of vertices in zone $z \in Z$, then $\mu_{q(v_i,v_j)}$ can be defined as follows:

$$\mu_{q(v_i,v_j)} = \max\{\mu_z, \mu_{z'}\}, v_i \in V_z, v_j \in V_{z'}, q(v_i, v_j) \in Q. \quad (27)$$

This approach prevents dilution of equity effects across mixed zones, ensures responsiveness to disadvantaged communities, and reflects fairness-oriented planning. It also offers transparency, as policymakers can trace each station's placement to the highest priority zone it serves.

*4.3. Mathematical Formulation*

The proposed model adopts a two-stage stochastic programming framework to optimize the placement of FCSs and MCSs under uncertainty. By explicitly incorporating spatial, temporal, and socioeconomic factors, the model promotes efficient use of financial resources and equitable access to EV charging infrastructure. The used notation and mathematical formulation of the first stage of the problem are presented next.

***Sets***

$Q$: Set of O-D pairs not covered by existing FCSs

$EP$: Set of endpoints

***Parameters***

$A$: $(|EP| \times |Q|)$ binary matrix, where element $a_{w\,q(v_i,v_j)}$ takes the value 1 if O-D pair $q(v_i, v_j) \in Q$ can be covered by a station located at endpoint $w \in EP$, and 0 otherwise

$c^f$: Cost of purchasing and installing a single FCS

$B^f$ : Budget allocated for purchasing and installing FCSs

$f(v_i, v_j)$: Average traffic flow between O-D pair $q(v_i, v_j)$ in roundtrips per time unit

$\mu_{q(v_i,v_j)}$: Composite weight for O-D pair $q(v_i, v_j) \in Q$

***First-stage decision variables***

$x_w$: Binary variable that takes the value 1 if an FCS is located at endpoint $w \in EP$, and 0 otherwise

$y_{q(v_i,v_j)}$: Binary variable that takes the value 1 if the O-D pair $q(v_i, v_j) \in Q$ is covered by at least one FCS, and 0 otherwise

The Fixed and Mobile Charging Location-Allocation problem can be formulated as a two-stage model. The first stage focuses on the placement of FCS to maximize a hybrid objective function maximizing coverage while considering socioeconomic weights in the network. The first stage mathematical model is as follows:

$$\text{Max: } F_1 = \sum_{q(v_i,v_j)\in Q} \left(1 + \mu_{q(v_i,v_j)}\right) f(v_i, v_j)\, y_{q(v_i,v_j)}, \quad (28)$$

Subject to:

$$\sum_{w\in EP} a_{wq(v_i,v_j)} x_w \geq y_{q(v_i,v_j)}, \forall q(v_i, v_j) \in Q, \quad (29)$$

$$\sum_{w\in EP} c^f x_w \leq B^f, \quad (30)$$

$$y_{q(v_i,v_j)}, x_w \in \{0, 1\}, \forall q(v_i, v_j) \in Q, w \in EP. \quad (31)$$

Equation ((28) defines the first-stage objective function, which maximizes the total weighted traffic flow covered by FCSs. Each O-D pair $q(v_i, v_j) \in Q$ is weighted by its average traffic flow $f(v_i, v_j)$ and a composite socioeconomic factor $\mu_{q(v_i,v_j)}$. Constraint set ((29) ensures that the flow between an O-D pair $q(v_i, v_j) \in Q$ is considered covered only if at least one endpoint $w \in EP$ capable of serving this pair is selected to install an FCS. Constraint ((30) imposes a budget limit on FCS deployment, restricting the total installation cost to the available budget $B^f$. Constraint ((31) defines the domain of the decision variables. The binary variables $x_w$, $w \in EP$ represent the first stage decisions. Specifically, if $x_w = 1$ then we recommend installing a new FCS at point $w \in EP$. Moreover, for each O-D pair $q(v_i, v_j) \in Q$ where $y_{q(v_i,v_j)} = 1$ the flow is covered by at least one FCS. These first-stage decisions establish the fixed infrastructure and its coverage, which will be used in the second stage as known parameters.

In the second stage, we consider a set of scenarios, denoted by $S$, representing uncertainty in traffic flow across a set of time periods, denoted by $T$. Let $N$ be the combined set of endpoints and existing charging station locations, defined as:

$$N = U \cup E. \quad (32)$$

Also, let $N^f$ be the set of locations with FCSs (existing and new) and $N^m$ be the set of remaining endpoints, not selected for installing new FCSs. These two sets can be defined as:

$$N^f = U \cup \{w \in EP | x_w = 1, w \in EP\}, \quad (33)$$

$$N^m = N \backslash N^f. \quad (34)$$

Moreover, let the set of uncovered O-D pairs by any FCS be denoted by $Q^m$. This set can be defined as:

$$Q^m = Q \backslash \left\{ q(v_i, v_j) | y_{q(v_i,v_j)} = 1, q(v_i, v_j) \in Q \right\}. \quad (35)$$

The second stage model allocates MCSs to points in $N^m$ to cover additional flows that are not covered by FCSs. For the second stage, we consider the following assumptions:

I. MCSs incur a relocation cost per mile when moving between any two points in the network.

II. After servicing a location in $N^m$, each MCS must travel to an FCS and remain there for at least one period to recharge before it can be dispatched to another service location.

The first assumption accounts for relocation costs, including fuel consumption and travel time, and discourages unnecessary movements of MCS across the network. The second assumption prevents MCSs from traveling directly between service locations without recharging, ensuring operational feasibility and reflecting realistic deployment practices.

The second-stage objective function incorporates both flow coverage and relocation cost. Since these two terms are expressed in different units, they cannot be directly aggregated. To address this, we introduce the parameter $\varepsilon$ which represents the combined economic and social benefit per unit of served flow [130]. This transformation converts flow coverage into a monetary equivalent, allowing both components to be expressed in consistent units and enabling the maximization of net system benefit. Below we present a table of additional notation for the second stage, followed by the mathematical model.

***Sets***

| | |
|---|---|
| $Q^m$: | Set of O-D pairs not covered by existing FCSs |
| $T$: | Set of discrete time periods |
| $S$: | Set of traffic flow scenarios |
| $M$: | Set of mobile charging stations |
| $N$: | Combined set of endpoints and existing charging station locations |
| $N^f$: | Set of FCS locations, including both existing and newly selected locations |
| $N^m$: | Set of endpoints available for MCS deployment |

***Parameters***

| | |
|---|---|
| $\varepsilon$: | Economic and social benefit per unit of served flow. |
| $\lambda$: | Cost of MCS relocation per mile |
| $P_s$: | Probability of scenario $s \in S$ |
| $f_{ts}(v_i, v_j)$: | Expected traffic flow between O-D pair $q(v_i, v_j) \in Q$ during period $t \in T$ under scenario $s \in S$ |
| $d(x_1, x_2)$: | Distance between points $x_1$ and $x_2$ in the network |
| $c^m$: | Cost of purchasing and installing a single MCS |
| $B^m$: | Budget allocated for purchasing MCSs |

***Second-stage decision variables***

| | |
|---|---|
| $\varphi_m$: | Binary variable that takes 1 if MCS $m$ is deployed at least once over the planning horizon, and 0 otherwise. |
| $k_{tmw}$: | Binary variable that takes the value 1 if MCS $m \in M$ is allocated to endpoint $w \in N$ in period $t \in T$, and 0 otherwise |
| $y_{tq(v_i,v_j)}$: | Binary variable that takes the value 1 if the O-D pair $q(v_i, v_j) \in Q$ is covered by at least one MCS in period $t \in T$, and 0 otherwise |
| $r_{tmww'}$: | Binary variable that takes the value 1 if MCS $m$ is relocated from endpoint $w$ in period $t$ to endpoint $w'$ in period $t+1$, where $w \in N$, and 0 otherwise |

The second stage mathematical model is as follows:

$$\text{Max: } F_2 = \sum_{s \in S} \sum_{t \in T} P_s \left( \sum_{q(v_i,v_j) \in Q^m} \varepsilon (1 + \mu_{q(v_i,v_j)}) f_{ts}(v_i, v_j) y_{tq(v_i,v_j)} - \sum_{m \in M} \sum_{w,w' \in N} \lambda r_{tmww'} d(w, w') \right), \tag{36}$$

Subject to:

$$\sum_{m \in M} \sum_{w \in N^m} a_{wq(v_i,v_j)} k_{tmw} \geq y_{tq(v_i,v_j)}, \forall q(v_i, v_j) \in Q^m, t \in T, \tag{37}$$

$$\sum_{w \in N} k_{tmw} \leq 1, \forall m \in M, \forall t \in T, \tag{38}$$

$$r_{tmww'} \geq k_{tmw} + k_{(t+1)mw'} - 1, \forall t \in \{1,2, \dots, |T| - 1\}, m \in M, w, w' \in N, \tag{39}$$

$$\sum_{w' \in N^f} k_{(t+1)mw'} \geq \sum_{w \in N^m} k_{tmw}, \forall t \in \{1,2, \dots, |T| - 1\}, m \in M, \tag{40}$$

$$\varphi_m \geq k_{tmw}, \forall m \in M, t \in T, w \in N^m, \tag{41}$$

$$\sum_{m \in M} c^m \varphi_m \leq B^m, \tag{42}$$

$$k_{tmw} \in \{0,1\}, \forall t \in T, m \in M, \forall w \in N^m, \tag{43}$$

$$y_{tq(v_i,v_j)}, r_{tmww'},\ \varphi_m \in \{0, 1\}, \forall t \in T, q(v_i, v_j) \in Q^m, m \in M, w, w' \in N. \tag{44}$$

Equation ((36) defines the second-stage objective function, which maximizes the expected benefit across all scenarios and time periods from extending coverage to flows that were not covered by FCSs. The first term captures the socioeconomic value of serving uncovered O-D pairs from $Q^m$ using MCSs. Each flow is weighted by its expected demand $f_{ts}(v_i, v_j)$, the equity factor $\left(1+\mu_{q(v_i,v_j)}\right)$, and the benefit parameter $\varepsilon$, which converts flow coverage into a monetary equivalent. The second term subtracts MCS relocation costs. Constraint set ((37) ensures that the flow between an O-D pair $q(v_i, v_j) \in Q$ is considered covered in period $t$ only if at least one MCS is stationed at an endpoint $w \in N^m$ capable of serving that pair. Constraint ((38) restricts each MCS to be assigned to at most one location per period. Constraint set ((39) establishes the logical relationship between MCS location decisions and relocation events. It activates the relocation variable when an MCS moves from endpoint $w \in N$ in period $t \in \{1,2, \dots, |T| - 1\}$ to endpoint $w' \in N$ in period $t + 1$. Constraint set ((40) enforces mandatory recharging behavior. It requires that if an MCS provides service at a service node $w \in N^m$ in period $t$, then it must relocate to a fixed charging station $w' \in N^f$ in period $t + 1$ to recharge. Constraint set ((41) defines the activation variable $\varphi_m$, which becomes one if MCS $m$ is deployed at least once over the planning horizon. Constraint ((42) imposes a budget limit on MCS deployment. Finally, constraints sets ((43) and ((44) define the domains decision variables.

## 5. Computational Experiments

This section illustrates and evaluates the performance of the proposed two-stage stochastic optimization model for the joint deployment of fixed and mobile EV charging stations. We examine multiple weighting scenarios that prioritize different socioeconomic factors to assess their impact on infrastructure placement and service coverage. The proposed framework is applied to the real transportation network of Bexar County, Texas. As illustrated in Figure 10, existing high-power public EV charging stations are concentrated in the central and northern regions of the county, while southern and peripheral areas remain underserved. This spatial imbalance highlights the need for strategic infrastructure planning to improve network accessibility, enhance spatial equity, and increase overall charging system efficiency. The remainder of this section is organized as follows. Subsection 5.1 describes the data generation process, including the construction of a realistic zone-based transportation network, identification of existing charging infrastructure, and incorporation of socioeconomic indicators. Subsection 5.2 presents the experimental setup, including solver configurations and model parameters. Subsection 5.3 reports and analyzes the results, focusing on stage-wise decisions, policy sensitivity, and the operational value of mobile charging flexibility.

### *5.1. Data Generation*

Computational experiments are conducted on a network based on Bexar County, Texas. Core inputs including geographic layout, zoning boundaries, inter-zone distances, and socioeconomic indicators are obtained from real-world data sources, while traffic flows are synthetically generated to represent a wide range of operating conditions. The zonal structure follows school district boundaries, as shown in Figure 10, and is constructed using median household income and education level data from the Bexar County GIS Hub and Northeast Independent School District (NEISD) datasets [131, 132]. Existing public DC fast charging stations are identified using the PlugShare platform. Only universally accessible public sites are included [133], forming the baseline infrastructure for evaluating the proposed deployment strategies. The resulting transportation network consists of 50 vertices and 70 edges. Vehicles' driving range is $R = 10$ miles. The portion of fuel in a vehicle's tank at the origin and destination vertices is $\alpha = 0.5$. Traffic flows $f(v_i, v_j)$ in round trips per day, for each O-D pair $q(v_i, v_j), i < j$, are generated randomly across the network. The distances between all existing and potential charging station locations were obtained from Google Maps data, measured in miles to reflect real-world travel conditions across the Bexar County network.

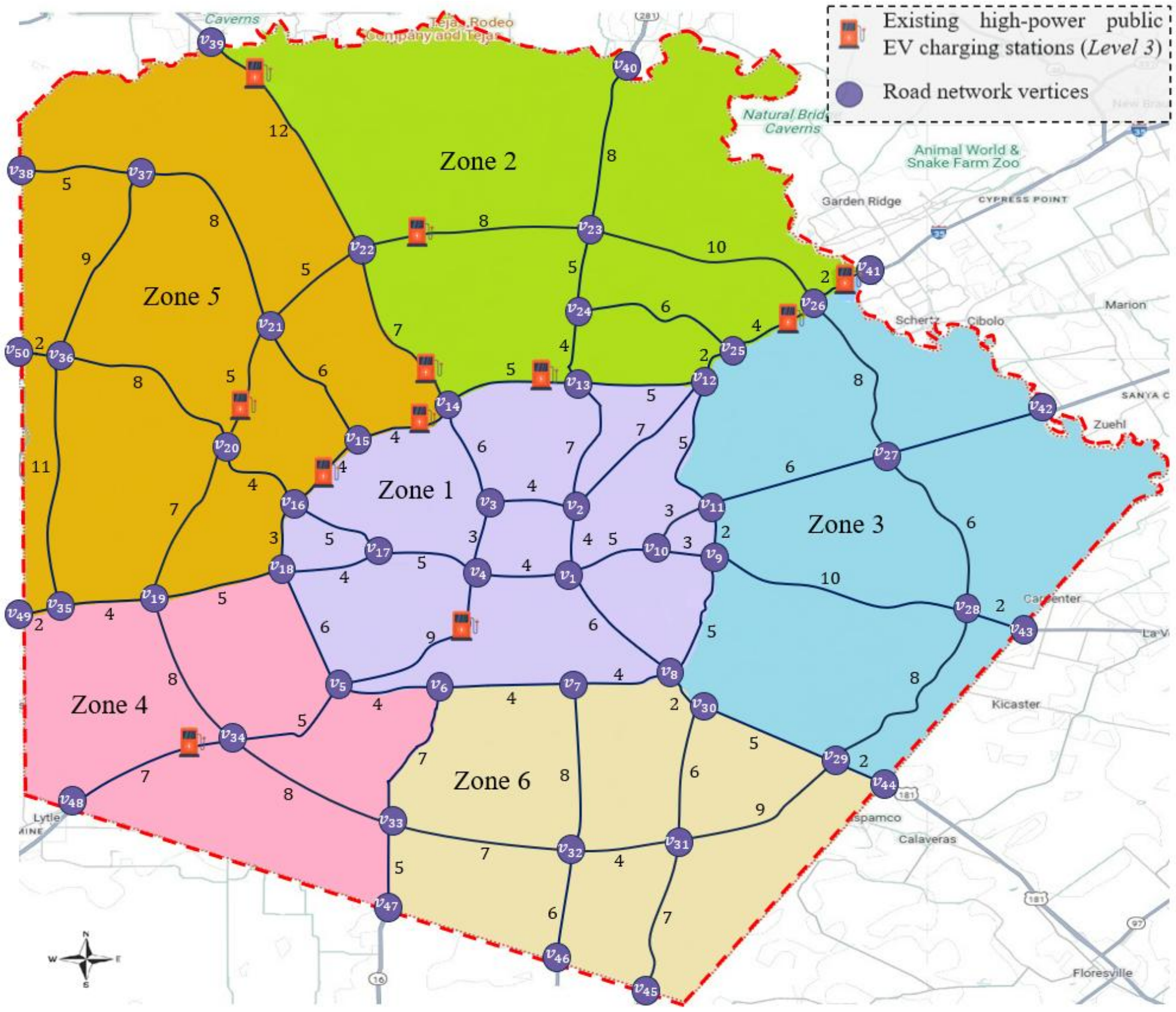

Figure 10: Distribution of High-Power Public EV Charging Stations Across Socioeconomic Zones in Bexar County

As discussed in previous sections, we consider four key socioeconomic factors to illustrate our proposed methodology, these factors are: median household income, population density, transit accessibility, and EV ownership rate. These factors capture both accessibility and adoption potential, and are widely employed in urban planning, public investment prioritization, and mobility justice frameworks [134].

Median household income and population density data are obtained from the U.S. Census Bureau's American Community Survey (ACS) and the City of San Antonio's open data portal [135-137]. Public transportation accessibility is assessed qualitatively based on proximity to major routes operated by VIA Metropolitan Transit [138]. Zones located near high-frequency and extensive bus corridors are classified as "High" accessibility, while areas with partial or limited service are categorized as "Medium" or "Low". EV ownership estimates are derived from vehicle registration records provided by the Texas Department of Motor Vehicles and aggregated trends reported by the Atlas EV Hub [139]. These data are combined with U.S. Census household statistics to compute an EV ownership rate, defined as the percentage of households within each zone that own at least one EV. Table 1 summarizes the values of the four socioeconomic factors across the six predefined zones in Bexar County.

Table 1: Socioeconomic Factors Across Defined Zones in Bexar County

| Zone | Median Income ($) | Population Density (per mi²) | Transit Access | EV Ownership (%) |
|---|---|---|---|---|
| $z$ | $\delta_1^z$ | $\delta_2^z$ | $\delta_3^z$ | $\delta_4^z$ |
| 1 | 42,000 | 5,200 | High | 4.5 |
| 2 | 60,000 | 3,800 | Medium | 5.1 |
| 3 | 37,000 | 6,000 | Low | 2.9 |
| 4 | 54,000 | 4,600 | Medium | 3.5 |
| 5 | 33,000 | 3,000 | Low | 2.2 |
| 6 | 48,000 | 5,800 | High | 3.8 |

*5.2. 5.2 Experimental Setup*

To evaluate the performance of the proposed optimization model, it was implemented in Python 3.11.10 and solved using Gurobi Optimizer 12.0.3. All computational experiments were conducted on a Windows workstation equipped with an Intel Core i9 processor (3.00 GHz) and 16 GB of RAM.

To examine different infrastructure planning strategies commonly used in practice by policymakers and experts, four cases are created using the BWM. In each case, experts first determine the most important (best) and least important (worst) socioeconomic factors based on policy goals and operational priorities. Using the standard 1-9 evaluation scale, expert judgments are then elicited to construct the Best-to-Others and Others-to-Worst comparison vectors $O_b$ and $O_w$, as defined in Equations ((18, (19). These assessments are then converted into quantitative factor weights using the BWM optimization model, with constraints set ((20-(24), which ensures consistency, transparency, and priority normalization.

In Case A (Equity-Focused Planning), policymakers emphasize economic vulnerability by selecting median household income as the most important factor and EV ownership as the least important, using $O_b = [1,3,6,9]$ and $O_w = [9,5,2,1]$. This configuration reflects a planning perspective that directs public investment toward economically disadvantaged communities. In Case B (Urban Access Prioritization), population density is identified as the primary driver of infrastructure need, while EV ownership remains the least important factor, with $O_b = [3,1,6,9]$ and $O_w = [5,9,2,1]$, capturing a strategy that prioritizes dense urban areas with high service demand. Case C (Mobility Justice Emphasis) shifts expert attention toward transit accessibility as the most important criterion, again ranking EV ownership as the least important, using $O_b = [6,3,1,9]$ and $O_w = [2,5,9,1]$. This case reflects policy objectives focused on supporting transit-dependent populations and improving equitable mobility outcomes. Finally, Case D (EV Adoption Readiness) represents a market-driven planning perspective in which EV ownership is treated as the most important factor and median income as the least, with $O_b = [9,3,6,1]$ and $O_w = [1,5,2,9]$, targeting zones with strong near-term adoption potential for rapid deployment. The resulting BWM-derived factor weights $\theta_i$, which reflect the relative importance assigned by policymakers and experts to each planning case, are summarized in Table 2.

Table 2: BWM-Derived Weights ($\theta_i$) for Socioeconomic Factors Across Four Expert Cases

| Case | Weight of Median Income | Weight of Population Density | Weight of Transit Access | Weight of EV Ownership |
|---|---|---|---|---|
| | $\theta_1$ | $\theta_2$ | $\theta_3$ | $\theta_4$ |
| A | 0.580 | 0.260 | 0.096 | 0.063 |
| B | 0.261 | 0.568 | 0.110 | 0.061 |
| C | 0.109 | 0.264 | 0.566 | 0.062 |
| D | 0.059 | 0.262 | 0.103 | 0.576 |

The raw socioeconomic values for each zone from Table 1 were normalized using Equation ((25) to ensure comparability across different scales. The resulting normalized values were then combined with the weight assignments $\theta_i$, derived from Table 2, to calculate the final zone weights $\mu_z$ using Equation ((26). These weights represent the relative social priority of each zone and are directly used in the proposed optimization framework. The normalized values of the socio-economic factors and the weights under each case are shown in Table 3. Next, we implement the Modified ESS algorithm to generate a set of candidate charging-station locations. In Step 1, we generate the set of O-D pairs $Q$. Note that some of the vertex pairs in Figure 10 are separated by distances greater than $R$; for example, $d(v_5, v_8) = 12$; therefore, the flow between this O-D pair cannot be covered by a single refueling station. At the same time, the flow $f(v_5, v_8) = 0$; hence, this O-D pair does not need to be considered. Using the definition of $Q$, the set of O-D pairs that need to be considered in this example consists of 104 O-D pairs, Appendix A (Table A.1 and Table A.2).

Table 3: Socioeconomic Factors Across Defined Zones in Bexar County

| Zone | Normalized Income | Normalized Density | Normalized Transit Score | Normalized EV Ownership | Case | Zone's Weight |
|---|---|---|---|---|---|---|
| | $(G_1^z)$ | $(G_2^z)$ | $(G_3^z)$ | $(G_4^z)$ | | $(\mu_z)$ |
| 1 | 0.67 | 0.73 | 0 | 0.79 | A | 0.63 |
| 2 | 0 | 0.27 | 0.5 | 1 | A | 0.18 |
| 3 | 0.85 | 1 | 1 | 0.24 | A | 0.86 |
| 4 | 0.22 | 0.53 | 0.5 | 0.45 | A | 0.34 |
| 5 | 1 | 0 | 1 | 0 | A | 0.67 |
| 6 | 0.44 | 0.93 | 0 | 0.55 | A | 0.53 |
| 1 | 0.67 | 0.73 | 0 | 0.79 | B | 0.64 |
| 2 | 0 | 0.26 | 0.5 | 1 | B | 0.26 |
| 3 | 0.86 | 1 | 1 | 0.24 | B | 0.92 |
| 4 | 0.23 | 0.53 | 0.5 | 0.44 | B | 0.44 |
| 5 | 1 | 0 | 1 | 0 | | 0.37 |
| 6 | 0.45 | 0.93 | 0 | 0.55 | | 0.68 |
| 1 | 0.67 | 0.73 | 0 | 0.79 | C | 0.31 |
| 2 | 0 | 0.26 | 0.5 | 1 | | 0.41 |
| 3 | 0.86 | 1 | 1 | 0.24 | | 0.94 |
| 4 | 0.23 | 0.53 | 0.5 | 0.44 | | 0.47 |
| 5 | 1 | 0 | 1 | 0 | | 0.65 |
| 6 | 0.45 | 0.93 | 0 | 0.55 | | 0.33 |
| 1 | 0.67 | 0.73 | 0 | 0.79 | D | 0.68 |
| 2 | 0 | 0.26 | 0.5 | 1 | | 0.70 |
| 3 | 0.86 | 1 | 1 | 0.24 | | 0.55 |
| 4 | 0.23 | 0.53 | 0.5 | 0.44 | | 0.46 |
| 5 | 1 | 0 | 1 | 0 | | 0.16 |
| 6 | 0.45 | 0.93 | 0 | 0.55 | | 0.59 |

In Step 2, for each edge $(v_a, v_b) \in E$, we consider each O-D pair $q(v_i, v_j) \in Q$ and the corresponding set of endpoints $EP(v_i, v_j; v_a, v_b)$ is found. To illustrate this process, let us consider edge $(v_{16}, v_{17})$ and O-D pair $q(v_{16}, v_{17})$; in this case $d(v_i, v_a) = d(v_{16}, v_{16}) = 0$ and $d(v_j, v_b) = d(v_{17}, v_{17}) = 0$. Using Equation ((6) for Type 1 case (a) trips:

$$\gamma(v_{16}; v_{16}, v_{17}) = \min\{5 - d(v_{16}, v_{16}), d(v_{16}, v_{17})\} = \min\{5 - 0,5\} = 5,$$

$$\gamma(v_{17}; v_{17}, v_{16}) = \min\{5 - d(v_{17}, v_{17}), d(v_{16}, v_{17})\} = \min\{5 - 4,5\} = 5.$$

Therefore, $RS_1^{16}(v_{16}, v_{17}; v_{16}, v_{17}) = \{x \in (v_{16}, v_{17}) | l(v_{16}, x) \leq 5 \text{ and } l(v_{17}, x) \leq 5\}$. The entire edge is a refueling segment for O-D pair $q(v_{16}, v_{17})$; therefore, there is no need to find other refueling segments for this O-D pair on edge $(v_{16}, v_{17})$.

Now, let us consider O-D pair $q(v_{16}, v_{18})$; in this case, $d(v_i, v_a) = d(v_{16}, v_{16}) = 0$ and $d(v_j, v_b) = d(v_{17}, v_{18}) = 4$. Using Equation (6) for Type 1 case (a) trips:

$$\gamma(v_{16}; v_{16}, v_{17}) = \min\{5 - d(v_{16}, v_{16}), d(v_{16}, v_{17})\} = \min\{5 - 0,5\} = 5,$$

$$\gamma(v_{18}; v_{17}, v_{16}) = \min\{5 - d(v_{18}, v_{17}), d(v_{16}, v_{17})\} = \min\{5 - 4,5\} = 1.$$

Therefore, $RS_1^{16}(v_{16}, v_{18}; v_{16}, v_{17}) = \{x \in (v_{16}, v_{17}) | l(v_{16}, x) \leq 5 \text{ and } l(v_{17}, x) \leq 1\}$. This refueling segment consists of a line segment along edge $(v_{16}, v_{17})$ starting from $v_{17}$ with length 1 and has two endpoints as shown in Figure 11.

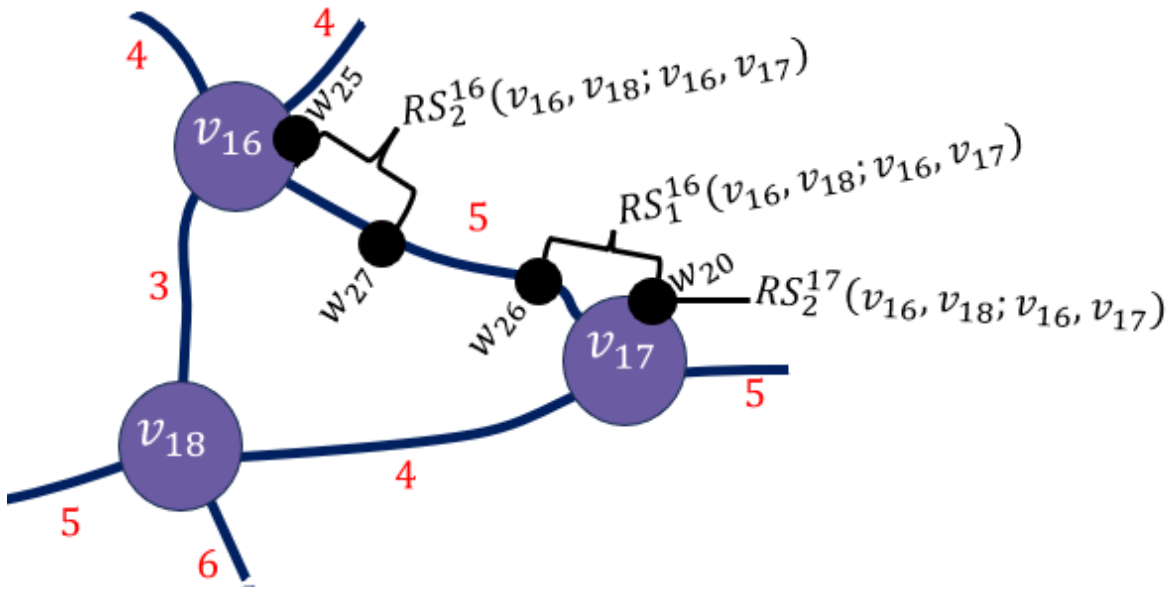

Figure 11: Refueling set of O-D pair $q(v_{16}, v_{18})$ on edge $(v_{16}, v_{17})$

Now, we find the refueling segment associated with Type 1 case (b) trip. Here, $d(v_i, v_b) = d(v_{16}, v_{17}) = 5$, $d(v_j, v_a) = d(v_{18}, v_{16}) = 3$, then use Equation (6) for Type 1 case (b) trips:

$$\gamma(v_{16}; v_{17}, v_{16}) = \min\{5 - d(v_{16}, v_{17}), d(v_{16}, v_{17})\} = \min\{5 - 5, 5\} = 0,$$

$$\gamma(v_{18}; v_{16}, v_{17}) = \min\{5 - d(v_{18}, v_{16}), d(v_{16}, v_{17})\} = \min\{5 - 3, 5\} = 2.$$

Therefore, $RS_1^{17}(v_{16}, v_{18}; v_{16}, v_{17}) = \{x \in (v_{16}, v_{17}) | l(v_{17}, x) \leq 0 \text{ and } l(v_{16}, x) \leq 2\} = \emptyset$. There is no refueling segment on edge $(v_{16}, v_{17})$ for Type 1 case (b) trips between O-D pair $q(v_{16}, v_{18})$. Next we use Equations ((9-(10)) to find the refueling segment for Type 2 case (a) trips as follows:

$$\delta^{16}(v_{16}, v_{18}; v_{16}, v_{17}) = \min\{5 - \max\{d(v_{16}, v_{16}), d(v_{18}, v_{16})\}, l(v_{16}, v_{17})\} = \min\{5 - \max\{0,3\},5\} = 2,$$

$$RS_2^{16}(v_{16}, v_{18}; v_{16}, v_{17}) = \{x \in (v_{16}, v_{17}) | l(v_{16}, x) \leq 2\},$$

Therefore, Type 2 case (a) refueling segment consists of a line segment with 2 units of length starting from $v_{16}$. Next, find the refueling segment for Type 2 case (b) trips:

$$\delta^{17}(v_{16}, v_{18}; v_{16}, v_{17}) = \min\{5 - \max\{d(v_{16}, v_{17}), d(v_{18}, v_{17})\}, l(v_{16}, v_{17})\} = \min\{5 - \max\{5,4\},5\} = 0,$$

$$RS_2^{17}(v_{16}, v_{18}; v_{16}, v_{17}) = \{x \in (v_{16}, v_{17}) \mid l(v_{17}, x) \leq 0\}.$$

Type 2 case (b) refueling segment consists of a single point at vertex $v_{17}$. Similarly, the remaining O-D pairs and edges in the network are analyzed. In Step 3, the sets of endpoints on all edges $(v_a, v_b) \in E$ are combined. The complete set of endpoints for this example, along with their locations on the network, is shown in Figure 12. Table  lists the potential charging station locations identified by the modified ESS algorithm within each zone of the Bexar County network.

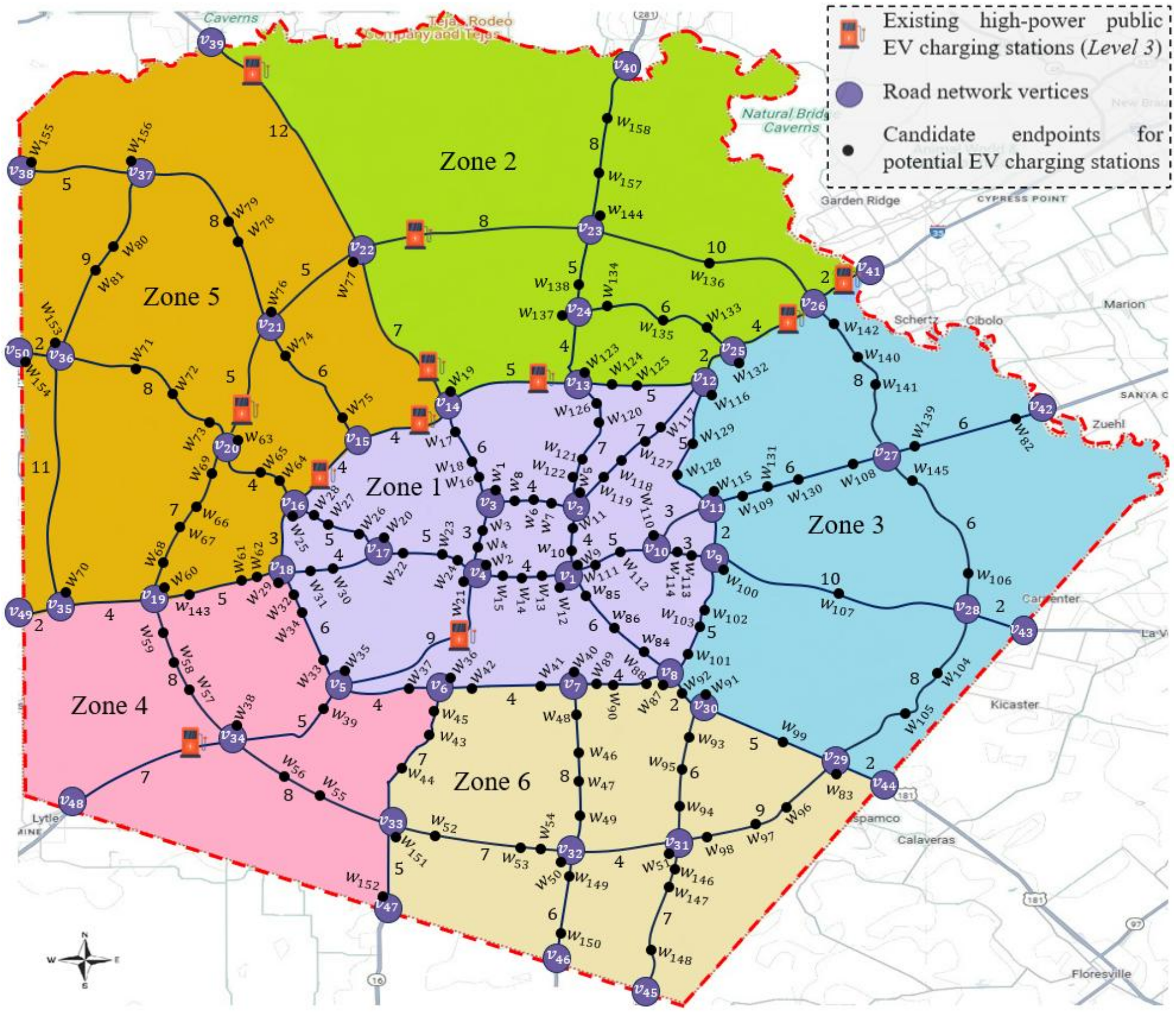

Figure 12: Complete set of endpoints Across Socioeconomic Zones in Bexar County

Table 4: Potential Charging Station Locations Within Each Zone in the Bexar County Network

| Zone | Set of candidate locations within zone $z$ |
|---|---|
| 1 | $w_1, w_2, w_3, w_4, w_5, w_6, w_7, w_8, w_9, w_{10}, w_{11}, w_{12}, w_{13}, w_{14}, w_{15}, w_{16}, w_{17}, w_{18}, w_{20}, w_{21}, w_{22}, w_{23}, w_{24}, w_{26}, w_{27}, w_{28}, w_{30}, w_{31}, w_{84}, w_{85}, w_{86}, w_{110}, w_{111}, w_{112}, w_{113}, w_{114}, w_{117}, w_{118}, w_{119}, w_{120}, w_{121}, w_{122}, w_{126}, w_{127}$ |
| 2 | $w_{19}, w_{77}, w_{116}, w_{123}, w_{124}, w_{125}, w_{132}, w_{133}, w_{134}, w_{135}, w_{136}, w_{137}, w_{138}, w_{144}, w_{157}, w_{158}$ |
| 3 | $w_{100}, w_{101}, w_{102}, w_{103}, w_{104}, w_{105}, w_{106}, w_{107}, w_{108}, w_{109}, w_{115}, w_{128}, w_{129}, w_{130}, w_{131}, w_{139}, w_{140}, w_{141}, w_{142}, w_{145}, w_{82}$ |
| 4 | $w_{29}, w_{32}, w_{33}, w_{34}, w_{35}, w_{36}, w_{37}, w_{38}, w_{39}, w_{43}, w_{44}, w_{45}, w_{55}, w_{56}, w_{57}, w_{58}, w_{59}, w_{60}, w_{61}, w_{62}, w_{70}, w_{143}, w_{151}, w_{152}$ |
| 5 | $w_{25}, w_{63}, w_{65}, w_{66}, w_{67}, w_{68}, w_{69}, w_{71}, w_{72}, w_{73}, w_{74}, w_{75}, w_{76}, w_{78}, w_{79}, w_{80}, w_{81}, w_{153}, w_{154}, w_{155}, w_{156}$ |
| 6 | $w_{40}, w_{41}, w_{42}, w_{46}, w_{47}, w_{48}, w_{49}, w_{50}, w_{51}, w_{52}, w_{53}, w_{54}, w_{87}, w_{88}, w_{89}, w_{90}, w_{91}, w_{92}, w_{93}, w_{94}, w_{95}, w_{96}, w_{97}, w_{98}, w_{99}, w_{146}, w_{147}, w_{148}, w_{149}, w_{150}, w_{83}$ |

### *5.3. Results and Analysis*

This section shows the results of the two-stage stochastic optimization model. The first stage distributes FCSs throughout Bexar County to maximize flow coverage while ensuring fair access through zone-based weights. The second stage adds MCSs to manage temporal demand changes and improve solution flexibility. Together, these results demonstrate the balance of efficiency, equity, and adaptability achieved under different weighting plans.

### *5.3.1. First Stage: FCS Allocation*

In the first stage, the optimization model identified the most appropriate endpoints for deploying FCSs under the four weighting Cases (A-D). The selected endpoints for each case are listed Table , and their spatial distributions are illustrated in Figure 13 and Appendix B (Figure B.1, Figure B.2, and Figure B.3), where blue EV charging station symbols indicate the selected endpoints for new FCSs.

Table 5: Selected Endpoints for FCS Allocation Across Four Weighting Cases (A–D)

| Case | Number of FCSs | Optimal Endpoint Locations | Allocated Zones | Flow Coverage Value |
|---|---|---|---|---|
| A | 10 | $w_{17}, w_{53}, w_{62}, w_{76}, w_{84}, w_{94}, w_{97}, w_{115}, w_{121}, w_{140}$ | $z_1, z_3, z_4, z_5, z_6$ | 582.8 |
| B | | $w_{32}, w_{40}, w_{52}, w_{77}, w_{84}, w_{93}, w_{97}, w_{115}, w_{121}, w_{141}$ | $z_1, z_2, z_3, z_4, z_6$ | 597.6 |
| C | | $w_{31}, w_{41}, w_{53}, w_{76}, w_{84}, w_{95}, w_{97}, w_{121}, w_{128}, w_{141}$ | $z_1, z_3, z_5, z_6$ | 584.1 |
| D | | $w_{18}, w_{32}, w_{52}, w_{77}, w_{84}, w_{95}, w_{97}, w_{115}, w_{121}, w_{141}$ | $z_1, z_2, z_3, z_4, z_6$ | 545.25 |

In Case A (Figure 13), the optimization model allocates endpoints across Zones 1, 3, 4, 5, and 6, with greater coverage in economically disadvantaged and densely populated areas. This spatial distribution reflects the equity-oriented weighting scheme. Although the total flow coverage (582.8) is lower than in other cases, the solution reduces spatial disparities across the network. Results and spatial distributions for Cases B, C, and D are provided in Appendix B.

### *5.3.2. Second Stage: MCS Deployment*

After establishing the FCSs, the second stage deploys MCSs to fill coverage gaps and respond to demand shifts. Travel flow is modeled dynamically across four periods and three demand scenarios: peak ($S_1, P_1 = 0.2$), shoulder ($S_2, P_2 = 0.5$), and off-peak ($S_3, P_3 = 0.2$), to reflect changes in charging demand. Under four weighting Cases (A-D), this stage assesses policy impacts on MCS placement. A deployment budget enough to have at most 10 MCSs is considered. Table presents the detailed spatiotemporal allocation and operational schedule of MCSs for the Equity-Focused Planning (Case A) scenario. In this case, MCSs are primarily deployed in low-income, high-demand zones such as Zones 1, 3, and 4, where permanent charging capacity is limited. The units cyclically alternate between serve and charge modes. No units are assigned to Zone 2, reflecting its relatively higher income and lower demand intensity. The results and operational schedules for the remaining three policy cases, Urban Access Prioritization (Case B), Mobility Justice Emphasis (Case C), and EV Adoption Readiness (Case D), are included in Appendix C, Table C.1. These cases support the findings of Case A by demonstrating the flexibility of the proposed framework under different strategic goals and demand scenarios, further confirming the model's robustness and policy significance.

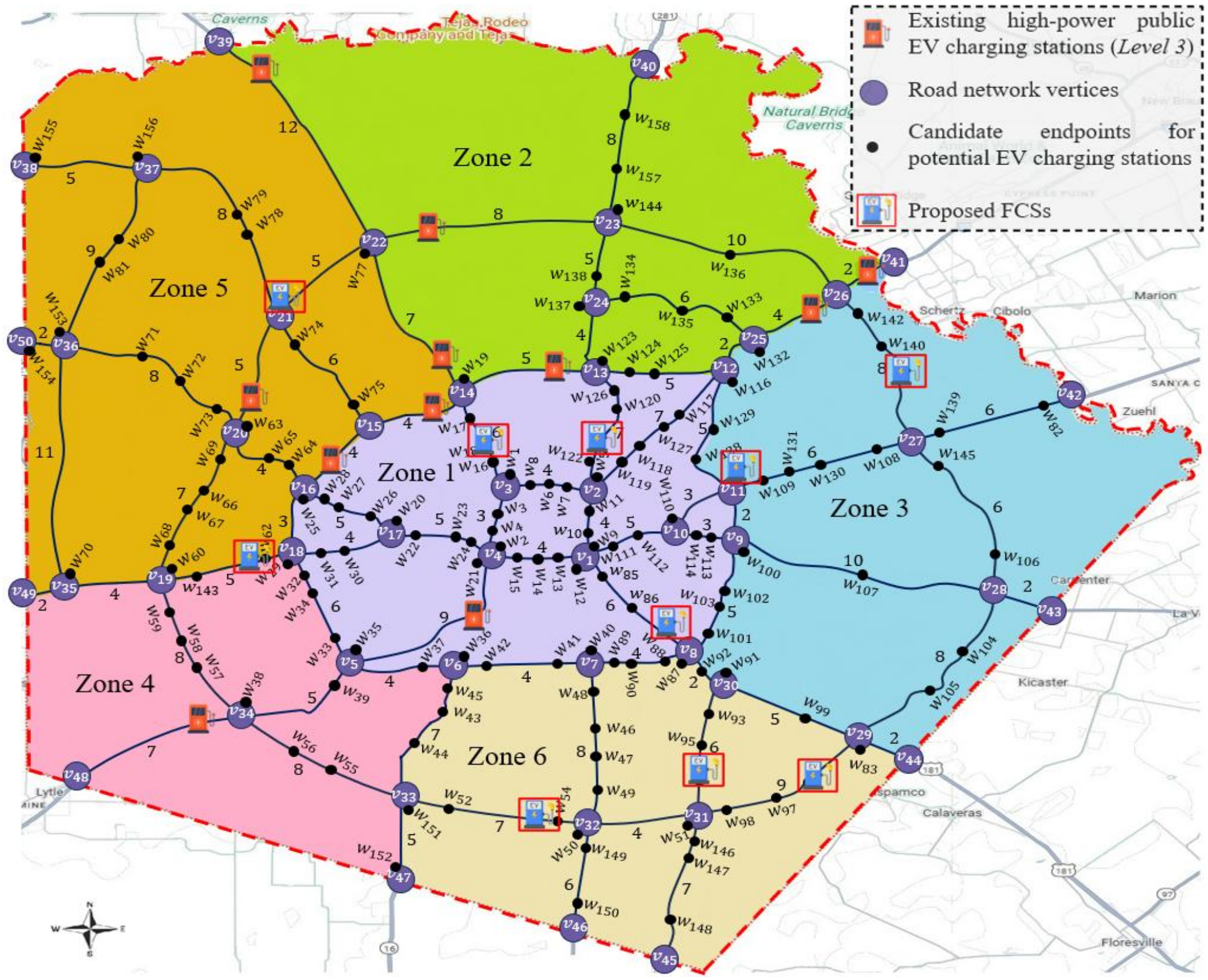

Figure 13: FCS allocation under Case A: Equity-Focused Planning

The coordinated serve–charge cycles observed in Case A illustrate the model's ability to manage MCS operations across multiple time periods under equity-oriented objectives. As shown in Figure 14, MCSs alternate between Serve and Charge modes across 6 zones and 4 time periods, with service deployments concentrated in higher-demand, lower-income areas such as Zones 1, 3, and 4. No MCS operates in consecutive service periods without recharging at a nearby FCS, ensuring energy feasibility while limiting unnecessary relocation. Zones with lower demand and higher socioeconomic status, such as Zone 2, receive no deployment. This spatiotemporal scheduling demonstrates how the optimization framework jointly coordinates service provision, recharging, and relocation to balance equity priorities, demand variation, and operational constraints.

Table 6: Spatiotemporal Allocation and Operational Schedules of MCSs Across Four Weighting Cases

| Case ID | Mobile Station ID | Time Period 1: Assigned Endpoint | Time Period 1: Allocated Zone | Time Period 1: Operation Mode | Time Period 2: Assigned Endpoint | Time Period 2: Allocated Zone | Time Period 2: Operation Mode | Time Period 3: Assigned Endpoint | Time Period 3: Allocated Zone | Time Period 3: Operation Mode | Time Period 4: Assigned Endpoint | Time Period 4: Allocated Zone | Time Period 4: Operation Mode |
|---|---|---|---|---|---|---|---|---|---|---|---|---|---|
| A | $M_1$ | $W_{84}$ | z1 | charge | $W_{12}$ | z1 | serve | $W_{84}$ | z1 | charge | $W_{87}$ | z6 | serve |
| | $M_2$ | $W_{29}$ | z4 | serve | $W_{62}$ | z4 | charge | $W_{29}$ | z4 | serve | $W_{62}$ | z4 | charge |
| | $M_3$ | $E_3$ | z3 | charge | $W_{75}$ | z5 | serve | $E_3$ | z3 | charge | $W_{25}$ | z5 | serve |
| | $M_4$ | $W_{87}$ | z6 | serve | $W_{84}$ | z1 | charge | $W_{87}$ | z6 | serve | $W_{84}$ | z1 | charge |
| | $M_5$ | $W_{53}$ | z6 | charge | $W_{149}$ | z6 | serve | $W_{53}$ | z6 | charge | $W_{149}$ | z6 | serve |
| | $M_6$ | $W_{62}$ | z4 | charge | $W_{29}$ | z4 | serve | $W_{62}$ | z4 | charge | $W_{60}$ | z4 | serve |
| | $M_7$ | $W_{84}$ | z1 | charge | $W_{87}$ | z6 | serve | $W_{84}$ | z1 | charge | $W_{99}$ | z6 | serve |
| | $M_8$ | $W_{110}$ | z1 | serve | $W_{115}$ | z3 | charge | $W_{109}$ | z3 | serve | $W_{115}$ | z3 | charge |
| | $M_9$ | $W_{115}$ | z3 | charge | $W_{109}$ | z3 | serve | $W_{115}$ | z3 | charge | $W_{109}$ | z3 | serve |
| | $M_{10}$ | $W_{25}$ | z5 | serve | $E_3$ | z3 | charge | $W_{25}$ | z5 | serve | $E_3$ | Z3 | charge |

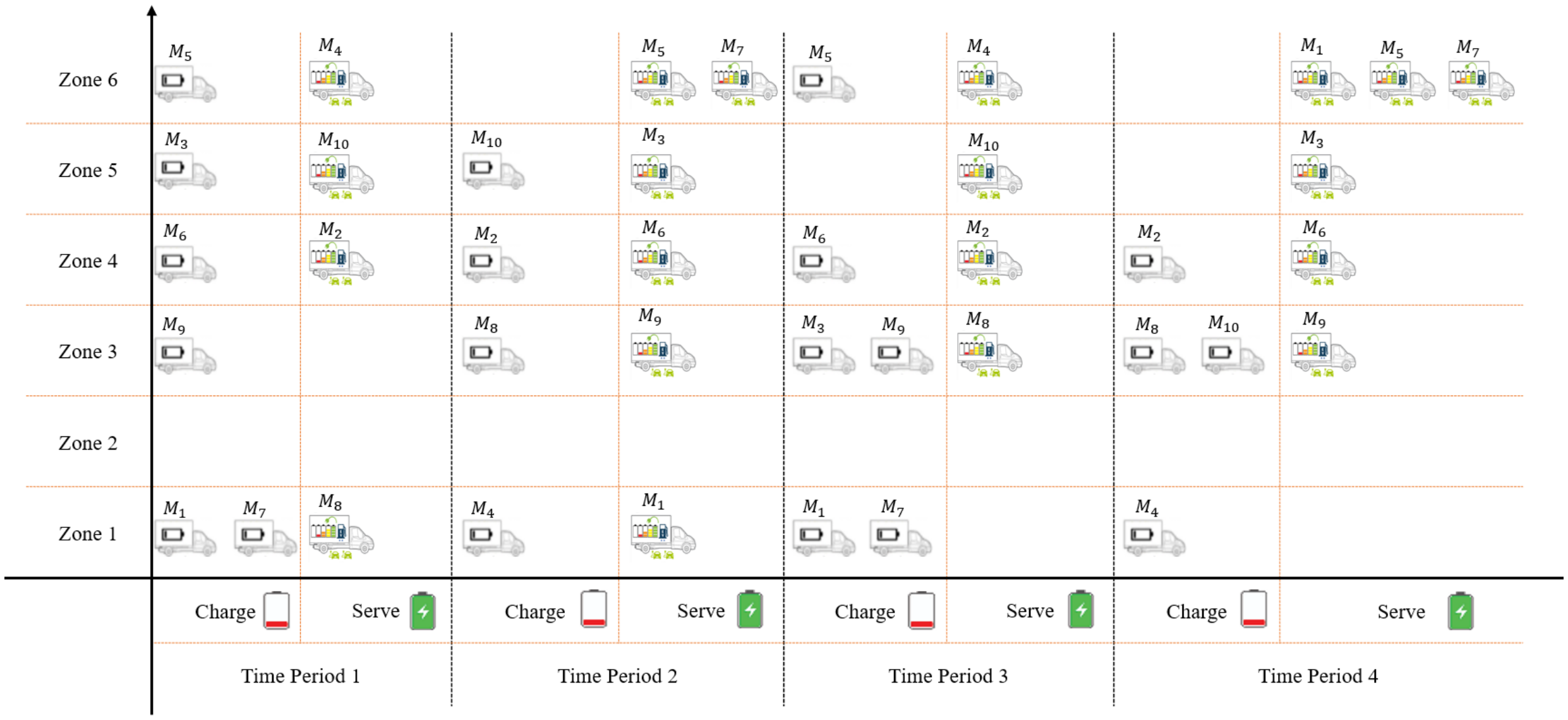

Figure 14: Temporal deployment of 10 MCSs across six zones during four time periods under Case A (Equity-Focused Planning)

Table  summarizes system performance across the four policy scenarios, incorporating both socio-economic benefits and relocation costs. Total benefits are computed using a socio-economic value of $\varepsilon$ = \$10 per unit of covered demand, capturing combined economic returns from EV charging and associated social benefits such as emission reductions and improved accessibility. Relocation costs are set at \$20 per

mile for moving an MCS between locations, representing the operational expense of repositioning units to maintain service continuity.

Table 7: System performance indicators across four different weighting Cases (A-D)

| Case ID | Objective Function Value | Total Socio-Economic Benefits | Relocation Cost |
|---|---|---|---|
| A | $ 8,008.62 | $ 9,208.63 | $ 1,200 |
| B | $ 4,149.70 | $ 7,089.70 | $ 2,940 |
| C | $ 7,203.88 | $ 8,843.88 | $ 1,640 |
| D | $ 3,659.34 | $ 6,019.34 | $ 2,360 |

Among the four cases, Case A (Equity-Focused Planning) yields the highest overall objective value ($8,008.62) and total benefit ($9,208.63), indicating that prioritizing lower-income, high-demand zones offers greater system-wide benefits with relatively low relocation costs ($1,200). Case C (Mobility Justice Emphasis) ranks second, striking a balance between social accessibility and operational efficiency. Next is Case B (Urban Access Prioritization), which incurs the highest relocation cost ($2,940) due to the frequent reallocation of MCSs in dense urban corridors. Case D (EV Adoption Readiness) offers the lowest benefit ($6,019.34) and objective value, as it focuses on high-income, high-ownership zones, which reduces the equity impact while maintaining moderate relocation costs ($2,360). Overall, these results indicate that socially driven policies (A and C) can achieve greater net system performance with fewer relocations, confirming the model's ability to align economic efficiency with equitable service delivery.

## 6. Conclusion

This paper presents a two-stage stochastic optimization framework for the equitable and resilient deployment of electric vehicle charging infrastructure through the coordinated use of fixed and mobile charging stations. The framework captures both spatial equity considerations and temporal demand variability by separating long-term planning decisions for FCSs from short-term operational decisions for MCS deployment. The first stage establishes a stable baseline network of FCSs using socioeconomic weights to prioritize underserved and high-demand areas, while a preprocessing step based on a modified ESS algorithm reduces redundancy and improves computational efficiency. The second stage enhances system adaptability by dynamically reallocating mobile charging stations across multiple periods, subject to operational feasibility and energy balance constraints.

Computational experiments on a realistic urban network demonstrate that different planning priorities lead to distinct trade-offs among coverage, relocation costs, and equity outcomes. Equity-focused and mobility-justice-oriented policies achieve higher net system benefits with fewer mobile relocations, indicating that socially driven planning strategies can improve overall performance without imposing excessive operational costs. In contrast, market-driven strategies that prioritize existing adoption potential yield lower equity impacts and reduced system-wide benefits. These results highlight the value of integrating socioeconomic considerations directly into infrastructure planning models and illustrate how mobile charging can complement fixed assets to improve resilience under uncertain and uneven demand.

Future research can extend this framework by incorporating variable charging prices to reflect time-dependent tariffs and user price sensitivity, as well as heterogeneous vehicle driving ranges to improve resource allocation accuracy. Modeling capacity limitations at both fixed and mobile charging stations would enable explicit representation of queueing effects and station utilization under high-demand conditions. Additional realism could be achieved by integrating grid-aware operational constraints, including power availability and coordination with energy providers. Together, these extensions would enhance the framework's ability to support policy analysis and long-term investment planning for large-scale EV charging infrastructure under uncertainty.

**Data availability**

The appendix presents the parameter values and distribution details for the randomly generated data used in the optimization model. For information on the simulation framework, please refer to the original citation.

## Appendix A.

Table A.1: Origin-Destination Pairs in the Network and Uncovered Demand Set.

| Set of all O-D pairs in the network | O-D pairs covered by existing stations |
|---|---|
| $q(v_1,v_2), q(v_1,v_4), q(v_1,v_8), q(v_1,v_{10}), q(v_2,v_3), q(v_2,v_{12}), q(v_2,v_{13}), q(v_3,v_4), q(v_3,v_{14}), q(v_4,v_5), q(v_4,v_{17}), q(v_5,v_6), q(v_5,v_{18}), q(v_5,v_{34}), q(v_6,v_7), q(v_6,v_{33}), q(v_7,v_8), q(v_7,v_{32}), q(v_8,v_9), q(v_8,v_{30}), q(v_9,v_{10}), q(v_9,v_{11}), q(v_9,v_{28}), q(v_{10},v_{11}), q(v_{11},v_{12}), q(v_{11},v_{27}), q(v_{12},v_{13}), q(v_{12},v_{25}), q(v_{13},v_{14}), q(v_{13},v_{24}), q(v_{14},v_{15}), q(v_{14},v_{22}), q(v_{15},v_{16}), q(v_{15},v_{21}), q(v_{16},v_{17}), q(v_{16},v_{18}), q(v_{16},v_{20}), q(v_{17},v_{18}), q(v_{19},v_{20}), q(v_{19},v_{34}), q(v_{19},v_{35}), q(v_{20},v_{21}), q(v_{20},v_{36}), q(v_{21},v_{22}), q(v_{21},v_{37}), q(v_{22},v_{23}), q(v_{22},v_{39}), q(v_{23},v_{24}), q(v_{23},v_{26}), (v_{23},v_{40}), q(v_{24},v_{25}), q(v_{25},v_{26}), q(v_{26},v_{27}), q(v_{26},v_{41}), q(v_{27},v_{28}), q(v_{27},v_{42}), q(v_{28},v_{29}), q(v_{28},v_{43}), q(v_{29},v_{30}), (v_{29},v_{31}), q(v_{29},v_{44}), q(v_{30},v_{31}), q(v_{31},v_{32}), q(v_{31},v_{45}), q(v_{32},v_{33}), q(v_{32},v_{46}), q(v_{33},v_{34}), q(v_{33},v_{47}), q(v_{34},v_{48}), q(v_{35},v_{36}), q(v_{35},v_{49}), q(v_{36},v_{37}), q(v_{36},v_{50}), q(v_{37},v_{38})$ | $q(v_4,v_5), q(v_{13},v_{14}), q(v_{14},v_{15}), q(v_{14},v_{22}), q(v_{15},v_{16}), q(v_{20},v_{21}), q(v_{22},v_{23}), q(v_{22},v_{39}), q(v_{25},v_{26}), q(v_{26},v_{41}), q(v_{34},v_{48})$ |

Table A.2: Subsets of All and Uncovered O-D Pairs for Each Zone z in the Bexar County Network.

| Zone | Set of all O-D pairs in the zone z | Set of uncovered O-D pairs in the zone z |
|---|---|---|
| 1 | $q(v_1,v_2), q(v_1,v_4), q(v_1,v_8), q(v_1,v_{10}), q(v_2,v_3), q(v_2,v_{12}), q(v_2,v_{13}), q(v_3,v_4), q(v_3,v_{14}), q(v_4,v_5), q(v_4,v_{17}), q(v_9,v_{10}), q(v_{10},v_{11}), q(v_{16},v_{17}), q(v_{17},v_{18})$ | $q(v_1,v_2), q(v_1,v_4), q(v_1,v_8), q(v_1,v_{10}), q(v_2,v_3), q(v_2,v_{12}), q(v_2,v_{13}), q(v_3,v_4), q(v_3,v_{14}), q(v_4,v_{17}), q(v_9,v_{10}), q(v_{10},v_{11}), q(v_{16},v_{17}), q(v_{17},v_{18})$ |
| 2 | $q(v_{12},v_{13}), q(v_{12},v_{25}), q(v_{13},v_{14}), q(v_{13},v_{24}), q(v_{14},v_{22}), q(v_{22},v_{23}), q(v_{22},v_{39}), q(v_{23},v_{24}), q(v_{23},v_{26}), q(v_{23},v_{40}), q(v_{24},v_{25}), q(v_{25},v_{26}), q(v_{26},v_{41})$ | $q(v_{12},v_{13}), q(v_{12},v_{25}), q(v_{13},v_{24}), q(v_{23},v_{24}), q(v_{23},v_{26}), q(v_{23},v_{40}), q(v_{24},v_{25})$ |
| 3 | $q(v_8,v_9), q(v_8,v_{30}), q(v_9,v_{11}), q(v_9,v_{28}), q(v_{11},v_{12}), q(v_{11},v_{27}), q(v_{26},v_{27}), q(v_{27},v_{28}), q(v_{27},v_{42}), q(v_{28},v_{29}), q(v_{28},v_{43}), q(v_{29},v_{30}), q(v_{29},v_{44})$ | $q(v_8,v_9), q(v_8,v_{30}), q(v_9,v_{11}), q(v_9,v_{28}), q(v_{11},v_{12}), q(v_{11},v_{27}), q(v_{26},v_{27}), q(v_{27},v_{28}), q(v_{27},v_{42}), q(v_{28},v_{29}), q(v_{28},v_{43}), q(v_{29},v_{30}), q(v_{29},v_{44})$ |
| 4 | $q(v_5,v_6), q(v_5,v_{18}), q(v_5,v_{34}), q(v_6,v_{33}), q(v_{18},v_{19}), q(v_{19},v_{34}), (v_{19},v_{35}), q(v_{33},v_{34}), q(v_{33},v_{47}), q(v_{34},v_{48}), q(v_{35},v_{49})$ | $q(v_5,v_6), q(v_5,v_{18}), q(v_5,v_{34}), q(v_6,v_{33}), q(v_{18},v_{19}), q(v_{19},v_{34}), q(v_{19},v_{35}), q(v_{33},v_{34}), q(v_{33},v_{47}), q(v_{35},v_{49})$ |
| 5 | $q(v_{14},v_{15}), q(v_{15},v_{16}), q(v_{15},v_{21}), q(v_{16},v_{18}), q(v_{16},v_{20}), q(v_{19},v_{20}), q(v_{20},v_{21}), q(v_{20},v_{36}), q(v_{21},v_{22}), q(v_{21},v_{37}), q(v_{35},v_{36}), q(v_{35},v_{49}), q(v_{36},v_{37}), q(v_{36},v_{50}), q(v_{37},v_{38})$ | $q(v_{15},v_{21}), q(v_{16},v_{18}), q(v_{16},v_{20}), q(v_{19},v_{20}), q(v_{20},v_{36}), q(v_{21},v_{22}), (v_{21},v_{37}), q(v_{35},v_{36}), q(v_{35},v_{49}), q(v_{36},v_{37}), q(v_{36},v_{50}), q(v_{37},v_{38})$ |
| 6 | $q(v_6,v_7), q(v_7,v_8), q(v_7,v_{32}), q(v_{29},v_{31}), q(v_{30},v_{31}), q(v_{31},v_{32}), q(v_{31},v_{45}), q(v_{32},v_{33}), q(v_{32},v_{46})$ | $q(v_6,v_7), q(v_7,v_8), q(v_7,v_{32}), q(v_{29},v_{31}), q(v_{30},v_{31}), q(v_{31},v_{32}), q(v_{31},v_{45}), q(v_{32},v_{33}), q(v_{32},v_{46})$ |

**Appendix B.** First-Stage Results: Spatial Distribution of FCSs Under Weighting Cases B-D. See Figure B.1, Figure B.2, and Figure B.3.

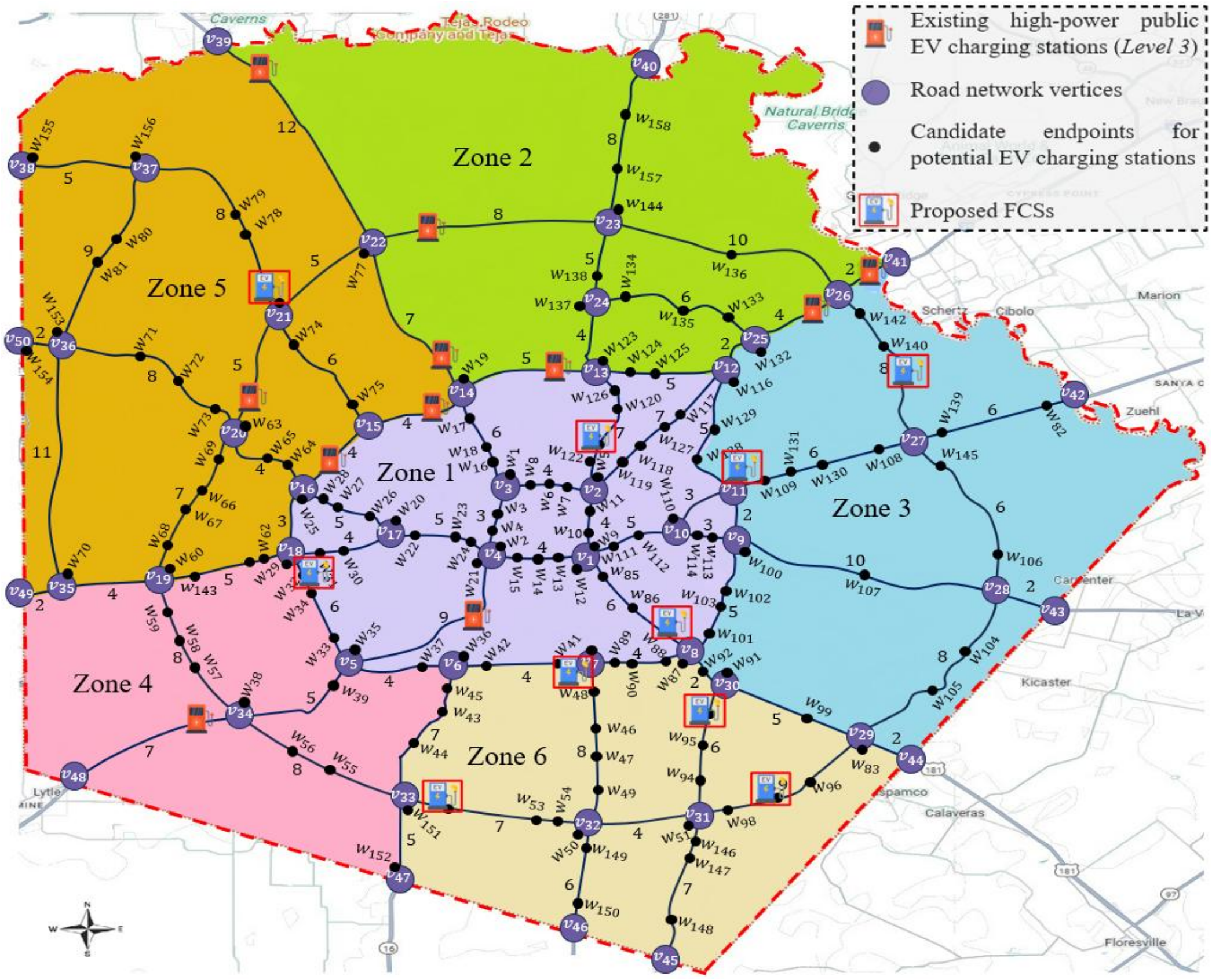

Figure B.1: FCS allocation under case B: Urban Access Prioritization

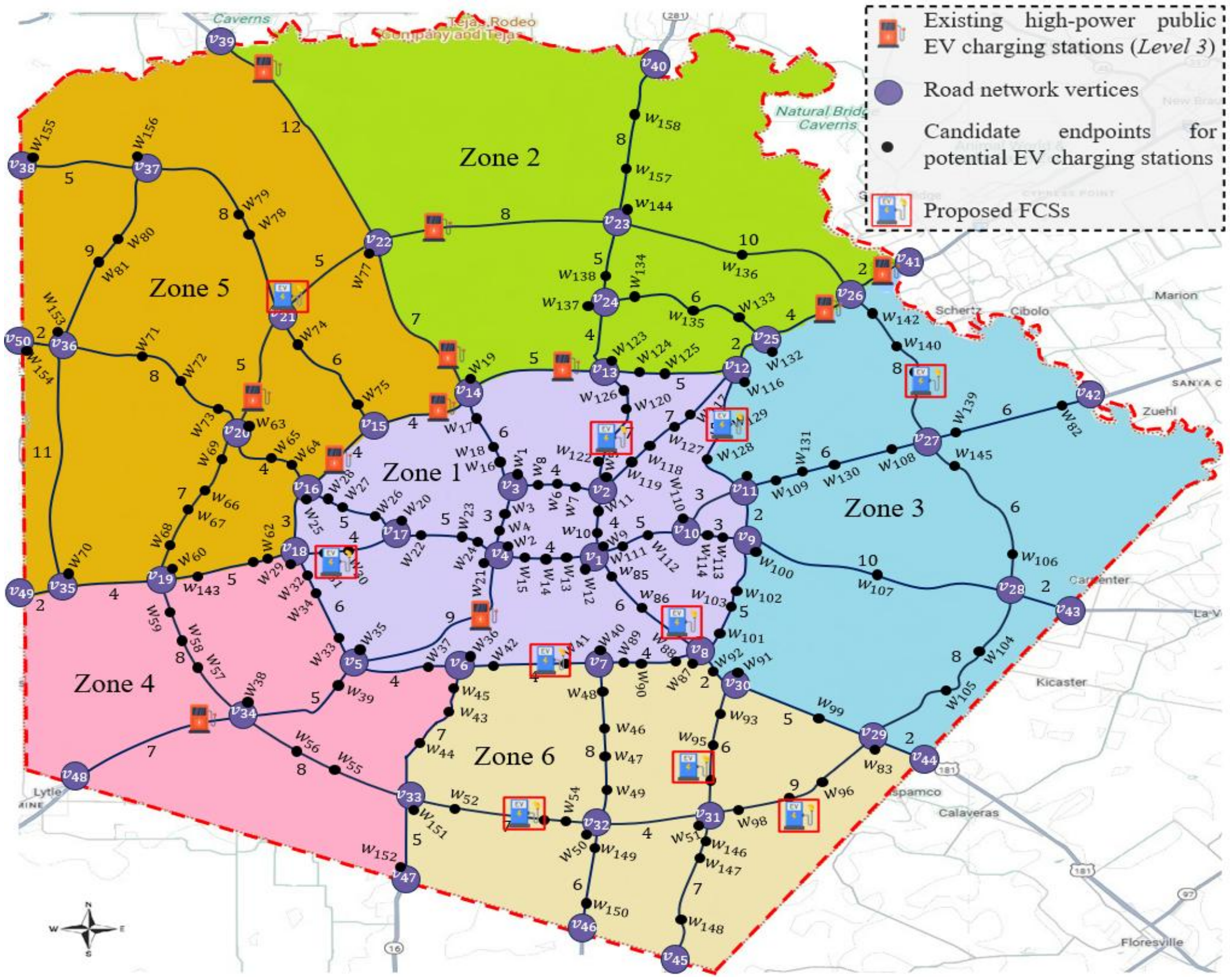

Figure B.2: FCS allocation under case C: Mobility Justice Emphasis

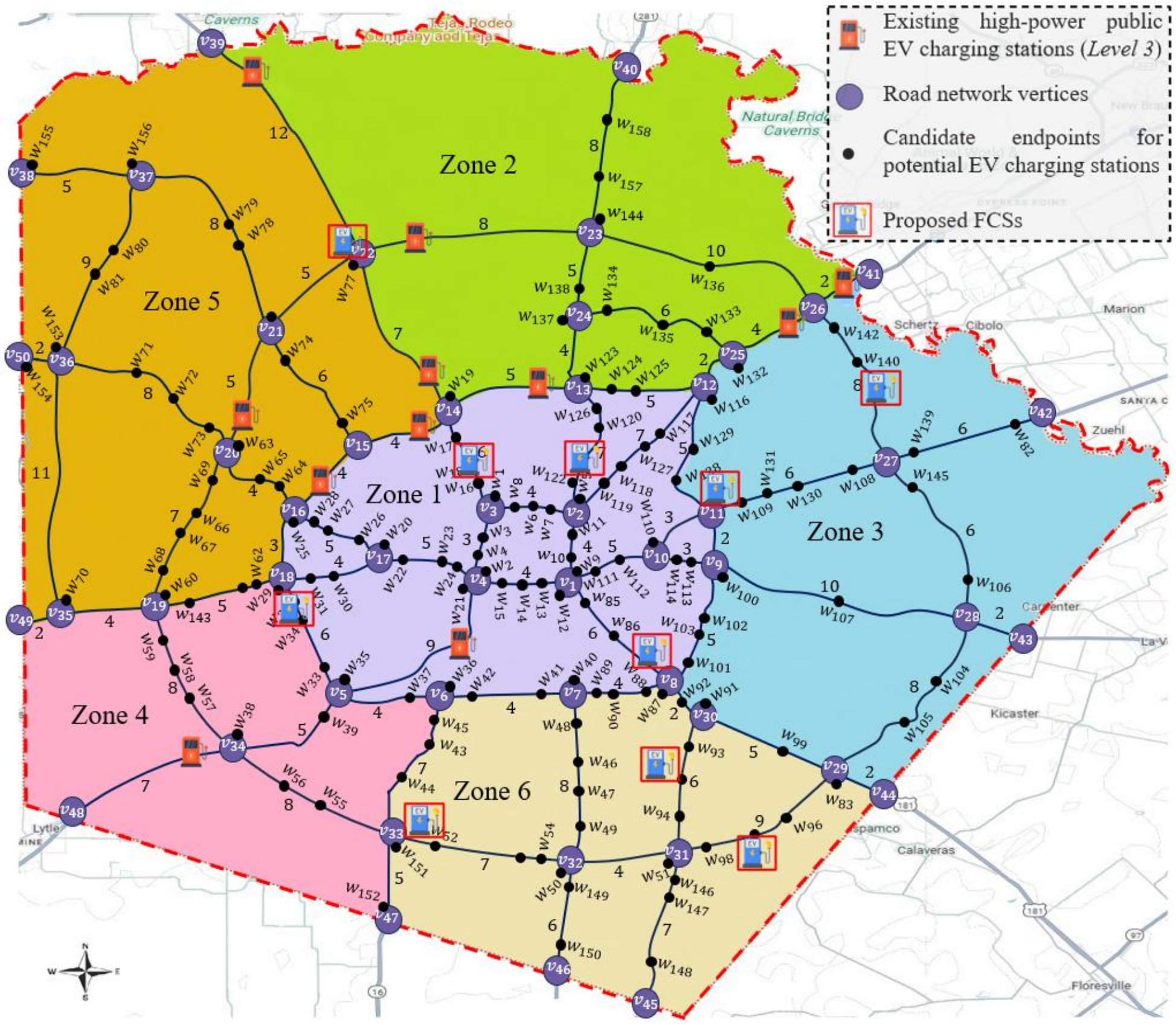

Figure B.3: FCS allocation under case D: EV Adoption Readiness

## Appendix C.

Table C.1: Spatiotemporal Allocation and Operational Schedules of MCSs for Weighting Cases (B-D)

| | | Time Period 1 | | | Time Period 2 | | | Time Period 3 | | | Time Period 4 | | |
|---|---|---|---|---|---|---|---|---|---|---|---|---|---|
| Case ID | Mobile Station ID | Assigned Endpoint | Allocated Zone | Operation Mode | Assigned Endpoint | Allocated Zone | Operation Mode | Assigned Endpoint | Allocated Zone | Operation Mode | Assigned Endpoint | Allocated Zone | Operation Mode |
| B | $M_1$ | $W_{104}$ | z3 | SERVE | $W_{97}$ | z6 | CHARGE | $W_{105}$ | z3 | SERVE | $W_{97}$ | z6 | CHARGE |
| | $M_2$ | $W_{149}$ | z6 | SERVE | $W_{52}$ | z6 | CHARGE | $W_{149}$ | z6 | SERVE | $W_{52}$ | z6 | CHARGE |
| | $M_3$ | $W_{121}$ | z1 | CHARGE | $W_5$ | z1 | SERVE | $W_{121}$ | z1 | CHARGE | $W_9$ | z1 | SERVE |
| | $M_4$ | $W_2$ | z1 | SERVE | $E_2$ | Z1 | CHARGE | $W_{35}$ | z4 | SERVE | $W_{32}$ | z4 | CHARGE |
| | $M_5$ | $W_{84}$ | z1 | CHARGE | $W_{87}$ | z6 | SERVE | $W_{84}$ | z1 | CHARGE | $W_{87}$ | z6 | SERVE |
| | $M_6$ | $E_2$ | z1 | CHARGE | $W_{21}$ | z1 | SERVE | $E_2$ | z1 | CHARGE | $W_{35}$ | z4 | SERVE |
| | $M_7$ | $W_5$ | z1 | SERVE | $W_{121}$ | z1 | CHARGE | $W_5$ | z1 | SERVE | $W_{121}$ | z1 | CHARGE |
| | $M_8$ | $W_{52}$ | z6 | CHARGE | $W_{149}$ | z6 | SERVE | $W_{52}$ | z6 | CHARGE | $W_{149}$ | z6 | SERVE |
| | $M_9$ | $W_{97}$ | z6 | CHARGE | $W_{105}$ | z3 | SERVE | $W_{97}$ | z3 | CHARGE | $W_{105}$ | z3 | SERVE |
| | $M_{10}$ | $W_{141}$ | z3 | CHARGE | $W_{139}$ | z3 | SERVE | $W_{141}$ | z5 | CHARGE | $W_{106}$ | z3 | SERVE |
| C | $M_1$ | $W_{84}$ | z1 | CHARGE | $W_{87}$ | z6 | SERVE | $W_{84}$ | z1 | CHARGE | $W_{87}$ | z6 | SERVE |
| | $M_2$ | $W_{128}$ | z3 | CHARGE | $W_{109}$ | z3 | SERVE | $W_{128}$ | z3 | CHARGE | $W_{106}$ | z3 | SERVE |
| | $M_3$ | $W_{87}$ | z6 | SERVE | $W_{84}$ | z1 | CHARGE | $W_{87}$ | z6 | SERVE | $W_{84}$ | z1 | CHARGE |
| | $M_4$ | $W_{41}$ | z6 | CHARGE | $W_{36}$ | z4 | SERVE | $W_{41}$ | z6 | CHARGE | $W_{36}$ | z4 | SERVE |
| | $M_5$ | $W_{109}$ | z3 | SERVE | $W_{128}$ | z3 | CHARGE | $W_{116}$ | z2 | SERVE | $W_{128}$ | z3 | CHARGE |
| | $M_6$ | $W_{97}$ | z6 | CHARGE | $W_{83}$ | z6 | SERVE | $W_{97}$ | z6 | CHARGE | $W_{83}$ | z6 | SERVE |
| | $M_7$ | $W_{29}$ | z4 | SERVE | $W_{31}$ | z1 | CHARGE | $W_{29}$ | z4 | SERVE | $W_{31}$ | z1 | CHARGE |
| | $M_8$ | $W_{116}$ | z2 | SERVE | $W_{128}$ | z3 | CHARGE | $W_{109}$ | z3 | SERVE | $W_{128}$ | z3 | CHARGE |
| | $M_9$ | $W_{128}$ | z3 | CHARGE | $W_{116}$ | z2 | SERVE | $W_{128}$ | z3 | CHARGE | $W_{116}$ | z2 | SERVE |
| | $M_{10}$ | $W_{31}$ | z1 | CHARGE | $W_{29}$ | z4 | SERVE | $W_{31}$ | z1 | CHARGE | $W_{29}$ | z4 | SERVE |
| D | $M_1$ | $E_{10}$ | z2 | CHARGE | $W_{133}$ | z2 | SERVE | $E_{10}$ | z2 | CHARGE | $W_{133}$ | z2 | SERVE |
| | $M_2$ | $W_{141}$ | z3 | CHARGE | $W_{139}$ | z3 | SERVE | $W_{141}$ | z3 | CHARGE | $W_{106}$ | z3 | SERVE |
| | $M_3$ | $E_2$ | z1 | CHARGE | $W_2$ | z1 | SERVE | $E_2$ | z1 | CHARGE | $W_2$ | z1 | SERVE |
| | $M_4$ | $W_{97}$ | z6 | CHARGE | $W_{83}$ | z6 | SERVE | $W_{97}$ | z6 | CHARGE | $W_{83}$ | z6 | SERVE |
| | $M_5$ | $W_{20}$ | z1 | SERVE | $W_{32}$ | z4 | CHARGE | $W_{35}$ | z4 | SERVE | $W_{32}$ | z4 | CHARGE |
| | $M_6$ | $W_2$ | z1 | SERVE | $E_2$ | z1 | CHARGE | $W_{21}$ | z1 | SERVE | $E_2$ | Z1 | CHARGE |
| | $M_7$ | $W_{149}$ | z6 | SERVE | $W_{52}$ | z6 | CHARGE | $W_{151}$ | z4 | SERVE | $W_{52}$ | z6 | CHARGE |
| | $M_8$ | $W_{52}$ | z6 | CHARGE | $W_{151}$ | z4 | SERVE | $W_{52}$ | z6 | CHARGE | $W_{149}$ | z6 | SERVE |
| | $M_9$ | $W_{32}$ | z4 | CHARGE | $W_{29}$ | z4 | SERVE | $W_{32}$ | z4 | CHARGE | $W_{35}$ | z4 | SERVE |
| | $M_{10}$ | $W_{83}$ | z6 | SERVE | $W_{97}$ | z6 | CHARGE | $W_{83}$ | z6 | SERVE | $W_{97}$ | z6 | CHARGE |